\documentclass[twocolumn]{autart}    

\usepackage{graphicx}          
                               
\usepackage{amsmath}
\usepackage{amssymb}

\newcommand*\diff{\mathop{}\!\mathrm{d}}%

\usepackage{graphics} 
\usepackage{epsfig} 
\usepackage{epstopdf}
\usepackage{subfigure}
\usepackage{color}

\graphicspath{{./figures/}}
\allowdisplaybreaks[4]

\begin{document}

\begin{frontmatter}

\title{ISS Property with Respect to Boundary Disturbances for a Class of Riesz-Spectral Boundary Control Systems\thanksref{footnoteinfo}} 

\thanks[footnoteinfo]{This publication has emanated from research supported in part by a research grant from Science Foundation Ireland (SFI) under grant number 16/RC/3872 and is co-funded under the European Regional Development Fund and by I-Form industry partners. 
The material in this paper was not presented at any conference.\\ 
Corresponding author H.~Lhachemi.}

\author[UCD]{Hugo Lhachemi}\ead{hugo.lhachemi@ucd.ie}, 
\author[UCD]{Robert Shorten}\ead{robert.shorten@ucd.ie},               

\address[UCD]{School of Electrical and Electronic Engineering, University College Dublin, Dublin, Ireland}  


\begin{abstract}                          
This paper deals with the establishment of Input-to-State Stability (ISS) estimates for infinite dimensional systems with respect to both boundary and distributed disturbances. First, a new approach is developed for the establishment of ISS estimates for a class of Riesz-spectral boundary control systems satisfying certain eigenvalue constraints. Second, a concept of weak solutions is introduced in order to relax the disturbances regularity assumptions required to ensure the existence of classical solutions. The proposed concept of weak solutions, that applies to a large class of boundary control systems which is not limited to the Riesz-spectral ones, provides a natural extension of the concept of both classical and mild solutions. Assuming that an ISS estimate holds true for classical solutions, we show the existence, the uniqueness, and the ISS property of the weak solutions.   
\end{abstract}

\end{frontmatter}

\section{Introduction}
The concept of Input-to-State Stability (ISS), originally introduced by Sontag for finite dimensional systems~\cite{sontag1989smooth}, is one of the main tools for assessing the robustness of a system with respect to external disturbances. This notion has been extensively investigated for finite dimensional nonlinear systems during the last three decades. More recently, the possible extension of ISS properties to both linear~\cite{argomedo2012d,karafyllis2016input,karafyllis2016iss,karafyllis2017iss,tanwani2017disturbance} and nonlinear~\cite{zheng2017giorgi,zheng2017input} Partial Differential Equations (PDEs), and more generally to infinite dimensional systems, has attracted much attention~\cite{karafyllis2019preview,mironchenko2016restatements,mironchenko2017characterizations}. 

For infinite dimensional systems, there exist essentially two distinct types of perturbations. The first type includes distributed perturbations; namely, perturbations acting directly in the state equation. The second type concerns boundary perturbations; namely, perturbations acting on the system state through an algebraic constraint by the means of an unbounded operator. In the case of PDEs, the distributed perturbations are also called in-domain perturbations as they appear directly in the PDEs. In contrast, the boundary perturbations appear in the boundary conditions of the PDEs. This second type of perturbation naturally appears in numerous boundary control problems such as heat equations~\cite{Curtain2012}, transport equations~\cite{karafyllis2016iss}, diffusion or diffusive equations~\cite{argomedo2012d}, and vibration of structures~\cite{Curtain2012} with numerous practical applications, e.g., in robotics~\cite{endo2017boundary,henikl2016infinite} and aerospace engineering~\cite{lhachemi2018boundaryAutomatica}.

While many results have been reported regarding the ISS property with respect to distributed disturbances~\cite{argomedo2013strict,dashkovskiy2011local,dashkovskiy2013input,mazenc2011strict,mironchenko2016local,mironchenko2014integral,mironchenko2015construction,prieur2012iss}, the establishment of ISS properties with respect to boundary disturbances remains challenging~\cite{argomedo2012d,karafyllis2016input,karafyllis2016iss,karafyllis2017iss}. 
The traditional method to study abstract boundary control systems consists of transfering the boundary disturbances into a distributed one by means of a bounded lifting operator. By doing so, the original boundary control system is made equivalent to a standard evolution equation for which efficient analysis tools exist. The main issue with such an approach is that the resulting distributed perturbation involves the time derivative of the boundary perturbation~\cite{Curtain2012}. In particular, this induces two main difficulties. First, this approach fails, in general, to establish the ISS property with respect to the boundary disturbance, but can only show the ISS property with respect to the first time derivative of the boundary disturbance. A notable exception is the case of monotone control systems for which the lifting approach can be successfully used to establish ISS estimates, see~\cite{mironchenko2019monotonicity} for details. Second, in order to ensure the existence of classical solutions, one has to assume that the boundary disturbance is twice continuously differentiable. The relaxation of this regularity assumption requires the introduction of a concept of mild or weak solutions extending the one of classical solutions. However, the explicit occurrence of the time derivative of the boundary perturbation in the evolution equation does not allow, in a general setting, a straightforward introduction of such a concept of mild or weak solutions for boundary disturbances that are only assumed, e.g., to be continuous~\cite{emirsjlow2000pdes}. 

Inspired by well established finite-dimensional techniques, it has been proposed to resort to Lyapunov functions to establish the ISS properties of PDEs~\cite{argomedo2012d,tanwani2017disturbance,zheng2017giorgi,zheng2017input}. An other approach, based on functional analysis tools, has been proposed in~\cite{karafyllis2016iss} for the study of 1-D parabolic equations. In the problem  therein, (the negative of) the underlying disturbance free operator belongs to the class of Sturm-Liouville operators. Thus, its eigenvectors can be selected such that they form a Hilbert basis of the underlying Hilbert space. By projecting the system trajectories onto this Hilbert basis and using the self-adjoint nature of the disturbance free operator, it was shown that the analysis of the system trajectories reduces to the study of a countably infinite number of Ordinary Differential Equations (ODEs). Each of these ODEs characterizes the time domain evolution of one coefficient of the system trajectory when expressed in the aforementioned Hilbert basis. Then, the ISS property was obtained by solving these ODEs and by resorting to Parseval’s identity. On one hand, this idea of projecting the system trajectories over a Hilbert basis was then investigated in~\cite{karafyllis2018boundary} for the study of the asymptotic gains of a damped wave equation. On the other hand, this approach was further investigated in~\cite{lhachemi2018input} for the study of the ISS properties of a clamped-free damped string. Speciffically, the system trajectories were projected over a Riesz basis~\cite{christensen2016introduction} formed by eigenvectors of the disturbance free operator. Then, the property connecting the norm of a vector and its coefficients in a Riesz basis was used to obtain the desired ISS estimates.

The contribution of this paper is twofold. First, we develop a new approach for the establishment of ISS estimates with respect to boundary disturbances for a class of Riesz-spectral boundary control systems satisfying certain eigenvalue constraints. The ISS property of such analytic semigroups was investigated first in~\cite{jacob2016input,jacob2018infinite,jacob2018continuity} for $L^\infty$ boundary inputs. The problem was embedded into the extrapolation space $\mathcal{H}_{-1}$ while invoking admissibility conditions for returning to the original Hilbert space $\mathcal{H}$. The approach adopted in this paper differs by generalizing the ideas developed first in~\cite{karafyllis2016iss} and then in~\cite{lhachemi2018input}. Assuming boundary and distributed disturbances of class $\mathcal{C}^2$ and $\mathcal{C}^1$, respectively, the ISS property is established for classical solutions by taking advantage of the projection of the system trajectories over a Riesz basis formed by the eigenvectors of the disturbance free operator. This new tool can be used either to derive general (see proof of Theorem~\ref{th: ISS classical solutions}) or system specific (see Subsection~\ref{subsec: Improvement in the neighborhood of alpha = 1}) versions of the ISS estimates. The advantage of this approach is fourfold. First, all the computations are performed within the original Hilbert space $\mathcal{H}$ while avoiding the embedding of the problem into the extrapolation space $\mathcal{H}_{-1}$. Second, it provides a spectral decomposition of the Riesz-spectral boundary control system under the form of a countable infinite number of ODEs describing the time domain evolution of the coefficients of the system trajectory in the Riesz basis. Such spectral decompositions have been used, e.g., for control design; see~\cite{coron2004global,coron2006global} for the feedback stabilization of heat and wave equations. However, these classical spectral decompositions involve the time derivative $\dot{d}$ of the boundary disturbance $d$. In sharp contrast, the one provided in this paper only involves the boundary perturbation $d$. Third, we show that the $d$-term of the aforementioned spectral decomposition can be related to stationary solutions of the studied Riesz-spectral boundary control system. Thus, the obtained constant related to the boundary perturbation in the ISS estimate is expressed as a function of the energy of stationary solutions of the abstract boundary control system. Finally, it is shown that the approach used for deriving the ISS estimate also allows the direct establishment of fading memory estimates without invoking a conversion lemma~\cite[Chap.~7]{karafyllis2019preview}. Such fading memory estimates are useful for the establishment of small gain conditions ensuring the stability of interconnected systems~\cite{karafyllis2019preview}.

The second contribution of this paper deals with the introduction of a concept of weak solutions for abstract boundary control systems that allows the relaxation of the perturbations regularity assumptions from $\mathcal{C}^2$ and $\mathcal{C}^1$ to $\mathcal{C}^0$. This approach applies to a general class of boundary control systems that is not limited to Riesz-spectral ones. The first concept of weak solutions for infinite dimensional nonhomogeneous Cauchy problems was origally introduced in~\cite{ball1977strongly}. In the case of Hilbert spaces, an alternative version of such a concept was proposed under a variational form, see, e.g., \cite[Def.~3.1.6]{Curtain2012}. The advantage of such a formulation is that it bridges the gap with the concept of solutions from distribution theory~\cite[A.5.29]{Curtain2012}. The purpose of this paper is to extend such a concept to the case of abstract boundary control systems. The proposed approach is an extension of the concept of classical solutions that presents three main advantages. First, the concept of weak solution is introduced direclty within the original Hilbert space $\mathcal{H}$ under a variational formulation. In particular, it provides an alternative framework that avoids either the traditional embedding of the original problem into the extrapolation space $\mathcal{H}_{-1}$~\cite{emirsjlow2000pdes,tucsnak2009observation} or the abstract extension of the mild solutions by density arguments~\cite{schmid2018stabilization}. Second, the proposed definition of weak solutions only depends on the operators of the abstract boundary control system; in particular, it does not make explicitly appear the $C_0$-semigroup generated by the disturbance free operator. Third, by its variational formulation and the introduction of test functions, the concept of weak solution is inspired by distribution theory. Assuming that an ISS estimate holds true with respect to the classical solutions of the abstract boundary control system, we show the existence and the uniqueness of the weak solutions, as well as their compatibility with the concept of mild solutions. It is also shown that the ISS estimate satisfied by classical solutions also holds true for weak solutions.

The remainder of this paper is organized as follows. Notations and definitions are introduced in Section~\ref{sec: notation, def, pb setting}. The establishment of an ISS estimate for a class of Riesz-spectral boundary control systems with respect to classical solutions is presented in Section~\ref{sec: ISS classical solutions}. Then, a concept of weak solutions and its properties are presented in Section~\ref{sec: ISS weak solutions} for a large class of boundary control systems. The obtained results are applied on an illustrative example in Section~\ref{sec: application}. Concluding remarks are provided in Section~\ref{sec: conclusion}.

\section{Notations and Definitions}\label{sec: notation, def, pb setting}

\subsection{Notations}
The sets of non-negative integers, positive integers, integers, real, non-negative real, positive real, negative real, and complex numbers are denoted by $\mathbb{N}$, $\mathbb{N}^*$, $\mathbb{Z}$, $\mathbb{R}$, $\mathbb{R}_+$, $\mathbb{R}_+^*$, $\mathbb{R}_-^*$, and $\mathbb{C}$, respectively. For any $z \in \mathbb{C}$, $\operatorname{Re} z$ denotes the real part of $z$. Throughout the paper, the field $\mathbb{K}$ is either $\mathbb{R}$ or $\mathbb{C}$. We consider the following classical classes of comparison function:
\begin{align*}
\mathcal{K} & = \left\{ \gamma \in \mathcal{C}^0(\mathbb{R}_+;\mathbb{R}_+) \,:\, \gamma \,\mathrm{strictly~increasing}, \, \gamma(0)=0 \right\}, \\
\mathcal{L} & = \left\{ \gamma \in \mathcal{C}^0(\mathbb{R}_+;\mathbb{R}_+) \,:\, \gamma \,\mathrm{strictly~decreasing}, \, \right. \\
& \phantom{=====================} \left. \underset{t \rightarrow +\infty}{\lim} \gamma(t) = 0 \right\}, \\
\mathcal{KL} & = \left\{ \beta \in \mathcal{C}^0(\mathbb{R}_+ \times \mathbb{R}_+;\mathbb{R}_+) \,:\, \beta(\cdot,t) \in \mathcal{K} \,\mathrm{for}\, t\geq 0, \, \right.\\
& \phantom{==================} \left. \beta(x,\cdot) \in \mathcal{L} \,\mathrm{for}\, x > 0 \right\}. 
\end{align*}
For an interval $I\subset\mathbb{R}$ and a $\mathbb{K}$-normed linear space $(E,\left\Vert \cdot \right\Vert_E)$, $\mathcal{C}^n(I;E)$ denotes the set of functions $f : I \rightarrow E$ that are $n$ times continuously differentiable. For any $a<b$, we endow $\mathcal{C}^0([a,b];E)$ with the norm $\Vert\cdot\Vert_{\mathcal{C}^0([a,b];E)}$ defined for any $f\in\mathcal{C}^0([a,b];E)$ by
\begin{equation*}
\Vert f \Vert_{\mathcal{C}^0([a,b];E)} \triangleq \underset{t \in [a,b]}{\sup} \Vert f(t) \Vert_E .
\end{equation*}
For a given linear operator $L$, $\mathrm{R}(L)$, $\mathrm{ker}(L)$, and $\rho(L)$ denote its range, its kernel, and its resolvent set, respectively. $\mathcal{L}(E,F)$ denotes the set of bounded linear operators from $E$ to $F$. Let $(\mathbb{K}^m,\Vert \cdot \Vert_{\mathbb{K}^m})$ be a normed space with $m \in \mathbb{N}^*$. For a given basis $\mathcal{E} = (e_1,e_2,\ldots,e_m)$ of $\mathbb{K}^m$, we denote by $\Vert \cdot \Vert_{\infty,\mathcal{E}}$ the infinity norm\footnote{Defined as the maximum of the absolute values of the vector components when written in the basis $\mathcal{E}$.} in $\mathcal{E}$. By virtue of the equivalence of the norms in finite dimension, we denote by $c(\mathcal{E})\in\mathbb{R}_+^*$ the smallest constant such that $\Vert\cdot\Vert_{\infty,\mathcal{E}} \leq c(\mathcal{E}) \Vert \cdot \Vert_{\mathbb{K}^m}$.

Finally, we introduce the Kronecker notation: $\delta_{a,b} = 1$ if $a=b$, $0$ otherwise. The time derivative of a real-valued differentiable function $f: I \rightarrow\mathbb{R}$ is denoted by $\dot{f}$. If $\mathcal{H}$ is a Hilbert space, the time derivative of an $\mathcal{H}$-valued differentiable function $f: I \rightarrow\mathcal{H}$ is denoted by $\mathrm{d} f / \mathrm{d}t$.

\subsection{Definitions and related properties}
In this paper, we consider the following definition of \emph{boundary control systems}~\cite[Def. 3.3.2]{Curtain2012}.

\begin{defn}[Boundary control system]
Let $(\mathcal{H},\left< \cdot , \cdot \right>_\mathcal{H})$ be a separable Hilbert space over $\mathbb{K}$ and $(\mathbb{K}^m,\Vert \cdot \Vert_{\mathbb{K}^m})$ be a $\mathbb{K}$-normed space with $m \in \mathbb{N}^*$. Consider the abstract system taking the form:
\begin{equation}\label{def: boundary control system}
\left\{\begin{split}
\dfrac{\mathrm{d} X}{\mathrm{d} t}(t) & = \mathcal{A} X(t) + U(t) , & t > 0 \\
\mathcal{B} X (t) & = d(t) , & t > 0 \\
X(0) & = X_0 
\end{split}\right.
\end{equation}
with
\begin{itemize}
\item $\mathcal{A} : D(\mathcal{A}) \subset \mathcal{H} \rightarrow \mathcal{H}$ a linear (unbounded) operator;
\item $\mathcal{B} : D(\mathcal{B}) \subset \mathcal{H} \rightarrow \mathbb{K}^m$ with $D(\mathcal{A}) \subset D(\mathcal{B})$ a linear boundary operator;
\item $U : \mathbb{R}_+ \rightarrow \mathcal{H}$ the distributed disturbance;
\item $d : \mathbb{R}_+ \rightarrow \mathbb{K}^m$ the lumped disturbance.
\end{itemize}
We say that $(\mathcal{A},\mathcal{B})$ is a boundary control system, with associated abstract system (\ref{def: boundary control system}), if
\begin{enumerate}
\item the disturbance free operator $\mathcal{A}_0$, defined over the domain $D(\mathcal{A}_0) \triangleq D(\mathcal{A}) \cap \mathrm{ker}(\mathcal{B})$ by $\mathcal{A}_0 \triangleq \left.\mathcal{A}\right|_{D(\mathcal{A}_0)}$, is the generator of a $C_0$-semigroup $S$ on $\mathcal{H}$;
\item there exists a bounded operator $B \in \mathcal{L}(\mathbb{K}^m,\mathcal{H})$, called a lifting operator, such that $\mathrm{R}(B) \subset D(\mathcal{A})$, $\mathcal{A}B \in \mathcal{L}(\mathbb{K}^m,\mathcal{H})$, and $\mathcal{B}B = I_{\mathbb{K}^m}$.
\end{enumerate}
\end{defn}

\begin{rem}
As the boundary input space $\mathbb{K}^m$ is finite dimensional, the bounded nature of $B$ and $\mathcal{A}B$ is immediate. Thus, the existence of the lifting operator $B$ reduces to the surjectivity of $\left.\mathcal{B}\right\vert_{D(\mathcal{A})}$.
\end{rem}

We then introduce the concept of Riesz-spectral operators~\cite[Def. 2.3.4]{Curtain2012}.

\begin{defn}[Riesz spectral operator]\label{def: Riesz-spectral operator}
Let $\mathcal{A}_0 : D(\mathcal{A}_0) \subset \mathcal{H} \rightarrow \mathcal{H}$ be a linear and closed operator with simple eigenvalues $\lambda_n$ and corresponding eigenvectors $\phi_n \in D(\mathcal{A}_0)$, $n \in \mathbb{N}$. $\mathcal{A}_0$ is a Riesz-spectral operator if
\begin{enumerate}
\item $\left\{ \phi_n , \; n \in \mathbb{N} \right\}$ is a Riesz basis~\cite{christensen2016introduction}:
\begin{enumerate}
\item $\left\{ \phi_n , \; n \in \mathbb{N} \right\}$ is maximal, i.e., $\overline{ \underset{n\in\mathbb{N}}{\mathrm{span}_\mathbb{K}} \;\phi_n } = \mathcal{H}$;
\item there exist constants $m_R, M_R \in \mathbb{R}_+^*$ such that for all $N \in \mathbb{N}$ and all $\alpha_0 , \ldots , \alpha_N \in \mathbb{K}$,
\begin{equation}\label{eq: Riesz basis - inequality}
m_R \sum\limits_{n=0}^{N} \vert \alpha_n \vert^2
\leq
\left\Vert \sum\limits_{n=0}^{N} \alpha_n \phi_n \right\Vert_\mathcal{H}^2
\leq
M_R \sum\limits_{n=0}^{N} \vert \alpha_n \vert^2 ;
\end{equation}
\end{enumerate}
\item the closure of $\{ \lambda_n , \; n \in \mathbb{N} \}$ is totally disconnected, i.e. for any distinct $a,b \in \overline{ \{ \lambda_n , \; n \in \mathbb{N} \} }$, $[a,b] \not\subset \overline{ \{ \lambda_n , \; n \in \mathbb{N} \} }$.
\end{enumerate}
\end{defn}

A subset of the properties satisfied by Riesz-spectral operators are gathered in the following lemma~\cite[Lemmas 2.3.2 and 2.3.3]{Curtain2012}.

\begin{lem}\label{lem: properties Riesz-spectral operators}
Let $\mathcal{A}_0 : D(\mathcal{A}_0) \subset \mathcal{H} \rightarrow \mathcal{H}$ be a Riesz-spectral operator. With the notations of Definition~\ref{def: Riesz-spectral operator}, the following are true.
\begin{itemize}
\item The eigenvalues of the adjoint operator $\mathcal{A}_0^*$ are provided for $n \in \mathbb{N}$ by $\mu_n \triangleq \overline{\lambda_n}$ and the associated eigenvectors $\psi_n \in D(\mathcal{A}_0^*)$ can be selected such that $\{ \phi_n , \; n\in\mathbb{N} \}$ and $\{ \psi_n , \; n\in\mathbb{N} \}$ are biorthogonal, i.e., for all $n,m\in\mathbb{N}$, $\left< \phi_n , \psi_m \right>_\mathcal{H} = \delta_{n,m}$.
\item The sequence of vectors $\{ \psi_n , \; n\in\mathbb{N} \}$ is a Riesz basis.
\item For all $(\alpha_n)_{n\in\mathbb{N}} \in \mathbb{K}^\mathbb{N}$, 
\begin{equation}\label{eq: riesz basis square summable sequence}
\sum\limits_{n\in\mathbb{N}} \vert \alpha_n \vert^2 < \infty
\Leftrightarrow
\sum\limits_{n\in\mathbb{N}} \alpha_n \phi_n \in \mathcal{H} . 
\end{equation}
\item For all $z \in \mathcal{H}$,
\begin{equation}\label{eq: riesz basis - decomposition}
z 
= \sum\limits_{n\in\mathbb{N}} \left< z , \psi_n \right>_\mathcal{H} \phi_n
= \sum\limits_{n\in\mathbb{N}} \left< z , \phi_n \right>_\mathcal{H} \psi_n .
\end{equation}
\item $\mathcal{A}_0$ is the generator of a $C_0$-semigroup $S$ if and only if $\underset{n \in \mathbb{N}}{\sup} \operatorname{Re} \lambda_n < \infty$. In this case, 
\begin{equation}\label{eq: Riesz-spectral operator - expression semigroup}
\forall t\in\mathbb{R}_+, \; \forall z \in \mathcal{H}, \;\; S(t)z = \sum\limits_{n\in\mathbb{N}} e^{\lambda_n t} \left< z , \psi_n \right>_\mathcal{H} \phi_n ,
\end{equation}
and its growth-bound satisfies $\omega_0 = \underset{n \in \mathbb{N}}{\sup} \operatorname{Re} \lambda_n$.
\end{itemize}
\end{lem}

Finally, we introduce the following definition.

\begin{defn}[Riesz-spectral boundary control system]
We say that $(\mathcal{A},\mathcal{B})$ is a Riesz-spectral boundary control system, with associated abstract system (\ref{def: boundary control system}), if
\begin{enumerate}
\item $(\mathcal{A},\mathcal{B})$ is a boundary control system with associated abstract system (\ref{def: boundary control system}); 
\item the underlying disturbance free operator $\mathcal{A}_0$ is a Riesz-spectral operator.
\end{enumerate} 
\end{defn}

In particular, we consider the class of Riesz-spectral boundary control systems such that: 
\begin{equation}\label{eq: problem setting - confinement conditions}
\omega_0 = \underset{n \in \mathbb{N}}{\sup} \; \operatorname{Re} \lambda_n < 0 
, \quad
\zeta \triangleq \underset{n \in \mathbb{N}}{\sup} \; \dfrac{\vert\lambda_n \vert}{\vert\operatorname{Re}\lambda_n \vert} < \infty .
\end{equation}
The constraint $\omega_0 < 0$ guarantees the exponential stability of the $C_0$-semigroup $S$ while $\zeta < +\infty$ ensures that the system is of parabolic type. The slight relaxion of the latter constraint is discussed in Remark~\ref{rem: theorem 1 relax zeta}.

\section{ISS Estimate for Classical Solutions}\label{sec: ISS classical solutions}

We investigate in this section the case of classical solutions~\cite[Def.~3.1.1]{Curtain2012} for (\ref{def: boundary control system}).  

\begin{defn}[Classical solutions]\label{def: classical solution}
Let $(\mathcal{A},\mathcal{B})$ be a boundary control system. Let $X_0 \in D(\mathcal{A})$, $d \in \mathcal{C}^0(\mathbb{R}_+;\mathbb{K}^m)$ such that $\mathcal{B}X_0 = d(0)$, and $U \in \mathcal{C}^0(\mathbb{R}_+;\mathcal{H})$ be given. We say that $X$ is a classical solution of (\ref{def: boundary control system}) associated with $(X_0,d,U)$ if $X \in \mathcal{C}^0(\mathbb{R}_+;D(\mathcal{A})) \cap \mathcal{C}^1(\mathbb{R}_+;\mathcal{H})$, $X(0)=X_0$, and for all $t \geq 0$, $(\mathrm{d}X/\mathrm{d}t)(t)=\mathcal{A}X(t)+U(t)$ and $\mathcal{B}X(t)=d(t)$. 
\end{defn}

\subsection{ISS for classical solutions}\label{subsec: ISS for classical solutions}

The ISS property of the studied class of Riesz-spectral boundary control systems was established in~\cite{jacob2016input,jacob2018infinite,jacob2018continuity} for $L^\infty$ boundary inputs by using the usual concept of mild solutions within the extrapolation space. In this Subsection~\ref{subsec: ISS for classical solutions}, we develop a different approach for assessing such a result for classical solutions. The proposed approach generalizes the ideas developed in~\cite{karafyllis2016iss,lhachemi2018input} consisting in the projection of the system trajectories over adequate Riesz bases. This approach relies on a novel spectral decomposition (see (\ref{eq: diff equation c_n}) and Remark~\ref{rem: novel spectral decomposition}) and can be used either to derive general (see proof of Theorem~\ref{th: ISS classical solutions}) or system specific (see Subsection~\ref{subsec: Improvement in the neighborhood of alpha = 1}) versions of the ISS estimates. Such an approach is further investigated in Subsection~\ref{subsec: energy-based ISS estimate} for the derivation of a novel energy-based ISS estimate (see Theorem~\ref{thm: constant C1 independant of B}), as well as the derivation of fading memory estimates (see Remark~\ref{rem: fading memory estimates}).

\begin{thm}\label{th: ISS classical solutions}
Let $(\mathcal{A},\mathcal{B})$ be a Riesz-spectral boundary control system such that the eigenvalue constraints (\ref{eq: problem setting - confinement conditions}) hold. For every initial condition $X_0 \in D(\mathcal{A})$, and every disturbance $d \in \mathcal{C}^2(\mathbb{R}_+;\mathbb{K}^m)$ and $U \in \mathcal{C}^1(\mathbb{R}_+;\mathcal{H})$ such that $\mathcal{B}X_0 = d(0)$, the abstract system (\ref{def: boundary control system}) has a unique classical solution $X \in \mathcal{C}^0(\mathbb{R}_+;D(\mathcal{A})) \cap \mathcal{C}^1(\mathbb{R}_+;\mathcal{H})$ associated with $(X_0,d,U)$. Furthermore, the system is exponentially ISS in the sense that there exist $C_0,C_1,C_2 \in \mathbb{R}_+^*$, independent of $X_0$, $d$, and $U$, such that for all $t \geq 0$, 
\begin{align}
\left\Vert X(t) \right\Vert_\mathcal{H}
\leq &
C_0 e^{-\kappa_0 t} \left\Vert X_0 \right\Vert_\mathcal{H} 
+ C_1 \Vert d \Vert_{\mathcal{C}^0([0,t],\mathbb{K}^m)} \nonumber \\
& + C_2 \Vert U \Vert_{\mathcal{C}^0([0,t];\mathcal{H})}  \label{eq: th classical solutions - ISS}
\end{align}
with $\kappa_0 = - \omega_0 > 0$, where $\omega_0$ is the growth bound of the $C_0$-semigroup $S$ generated by the disturbance free operator $\mathcal{A}_0$. 
\end{thm}

\textbf{Proof of Theorem~\ref{th: ISS classical solutions}.}
The existence and uniqueness of the classical solutions under the assumed regularity assumptions directly follows from classical results for abstract boundary control systems (see, e.g., \cite[Thm~3.1.3 and Thm~3.3.3]{Curtain2012}). Thus, the proof is devoted to the derivation of the ISS estimate (\ref{eq: th classical solutions - ISS}). Let $B$ be a lifting operator associated with the boundary control system $(\mathcal{A},\mathcal{B})$ as provided by Definition~\ref{def: boundary control system}. Let $X_0 \in D(\mathcal{A})$, $d \in \mathcal{C}^2(\mathbb{R}_+;\mathbb{K}^m)$, and $U \in \mathcal{C}^1(\mathbb{R}_+;\mathcal{H})$ such that $\mathcal{B}X_0 = d(0)$ be given. Let $X$ be the classical solution of (\ref{def: boundary control system}) associated with $(X_0,d,u)$. Adopting the notations of Definition~\ref{def: Riesz-spectral operator} and Lemma~\ref{lem: properties Riesz-spectral operators}, we have from (\ref{eq: riesz basis - decomposition}) that for all $t \geq 0$,
\begin{equation}\label{eq: projection of X into the Riesz basis}
X(t) = \sum\limits_{n \in \mathbb{N}} \left< X(t) , \psi_n \right>_\mathcal{H} \phi_n .
\end{equation}
Introducing for all $n \in \mathbb{N}$ and all $t \geq 0$, $c_n(t) \triangleq \left< X(t) , \psi_n \right>_\mathcal{H}$, then $c_n \in \mathcal{C}^1(\mathbb{R}_+;\mathbb{K})$ and its time derivative is given by
\begin{align*}
\dot{c}_n(t) 
& = \left< \dfrac{\mathrm{d} X}{\mathrm{d} t}(t) , \psi_n \right>_\mathcal{H} \\
& = \left< \mathcal{A} X(t) + U(t) , \psi_n \right>_\mathcal{H} \\
& = \left< \mathcal{A} \left\{ X(t) - B d(t) \right\} , \psi_n \right>_\mathcal{H} + \left< \mathcal{A} B d(t) , \psi_n \right>_\mathcal{H} \\
& \phantom{=}\; + \left< U(t) , \psi_n \right>_\mathcal{H} ,
\end{align*}
where the last equality holds true because $X(t) \in D(\mathcal{A})$ and $\mathrm{R}(B) \subset{D(\mathcal{A})}$, providing $X(t)-Bd(t) \in D(\mathcal{A})$. Furthermore, $\mathcal{B} \{ X(t) - Bd(t) \} = d(t) - d(t) = 0$. Thus $X(t)-Bd(t) \in D(\mathcal{A}_0)$, yielding 
\begin{align*}
\left< \mathcal{A} \left\{ X(t) - B d(t) \right\} , \psi_n \right>_\mathcal{H} 
& = \left< \mathcal{A}_0 \left\{ X(t) - B d(t) \right\} , \psi_n \right>_\mathcal{H} \\
& = \left< X(t) - B d(t) , \mathcal{A}_0^* \psi_n \right>_\mathcal{H} \\
& = \left< X(t) - B d(t) , \overline{\lambda_n} \psi_n \right>_\mathcal{H} \\
& = \lambda_n \left< X(t) - B d(t) , \psi_n \right>_\mathcal{H} . 
\end{align*}
We get for all $t \geq 0$,
\begin{align}
\dot{c}_n(t) 
& = \lambda_n c_n(t) - \lambda_n \left< B d(t) , \psi_n \right>_\mathcal{H} \nonumber \\
& \phantom{=}\; + \left< \mathcal{A} B d(t) , \psi_n \right>_\mathcal{H} + \left< U(t) , \psi_n \right>_\mathcal{H} . \label{eq: diff equation c_n}
\end{align}
As all the terms involved in (\ref{eq: diff equation c_n}) are continuous over $\mathbb{R}_+$, a straightforward integration gives for all $t \geq 0$,
\begin{align}
c_n(t) 
& = e^{\lambda_n t} c_n(0) - \lambda_n \int_0^t e^{\lambda_n (t-\tau)} \left< B d(\tau) , \psi_n \right>_\mathcal{H} \diff\tau \nonumber \\
& \phantom{=}\; + \int_0^t e^{\lambda_n (t-\tau)} \left< \mathcal{A} B d(\tau) , \psi_n \right>_\mathcal{H} \diff\tau \label{eq: integral expression c_n} \\
& \phantom{=}\; + \int_0^t e^{\lambda_n (t-\tau)} \left< U(\tau) , \psi_n \right>_\mathcal{H} \diff\tau . \nonumber
\end{align}
Note that
\begin{equation*}
X(t) = \sum\limits_{n \in \mathbb{N}} c_n(t) \phi_n , \qquad
S(t)X_0 = \sum\limits_{n \in \mathbb{N}} e^{\lambda_n t} c_n(0) \phi_n ,
\end{equation*}
\begin{align*}
& \int_0^t S(t-\tau) \mathcal{A} B d(\tau) \diff\tau \\
& \qquad = \sum\limits_{n \in \mathbb{N}} \left< \int_0^t S(t-\tau) \mathcal{A} B d(\tau) \diff\tau, \psi_n \right>_\mathcal{H} \phi_n \\
& \qquad = \sum\limits_{n \in \mathbb{N}} \int_0^t \left< S(t-\tau) \mathcal{A} B d(\tau) , \psi_n \right>_\mathcal{H} \diff\tau \,\phi_n \\
& \qquad = \sum\limits_{n \in \mathbb{N}} \int_0^t e^{\lambda_n (t-\tau)} \left< \mathcal{A} B d(\tau) , \psi_n \right>_\mathcal{H} \diff\tau \,\phi_n ,
\end{align*}
and, similarly,
\begin{equation}\label{eq: trajectory - integral term U}
\int_0^t S(t-\tau) U(\tau) \diff\tau  = \sum\limits_{n \in \mathbb{N}} \int_0^t e^{\lambda_n (t-\tau)} \left< U(\tau) , \psi_n \right>_\mathcal{H} \diff\tau \,\phi_n .
\end{equation}
Thus, introducing
\begin{equation}\label{eq: coeff alpha(t)}
\alpha_n(t) \triangleq \lambda_n \int_0^t e^{\lambda_n (t-\tau)} \left< B d(\tau) , \psi_n \right>_\mathcal{H} \diff\tau , 
\end{equation}
we deduce from (\ref{eq: riesz basis square summable sequence}) and (\ref{eq: integral expression c_n}) that $(\alpha_n(t))_{n \in \mathbb{N}}$ is a square summable sequence for all $t \geq 0$ and that
\begin{equation*}
\alpha(t) \triangleq \sum\limits_{n \in \mathbb{N}} \alpha_n(t) \phi_n \in \mathcal{H} . 
\end{equation*}
Therefore, multiplying both sides of (\ref{eq: integral expression c_n}) by $\phi_n$ and summing over $n \in \mathbb{N}$ yields
\begin{align*}
X(t) & =  S(t)X_0 - \alpha(t) \nonumber \\
& \phantom{=}\; + \int_0^t S(t-\tau) \left\{  \mathcal{A} B d(\tau) + U(\tau) \right\} \diff\tau .
\end{align*}
Thus, for all $t \geq 0$, we have 
\begin{align}
\Vert X(t) \Vert_\mathcal{H} 
& \leq \Vert S(t) X_0 \Vert_\mathcal{H} 
+ \Vert \alpha(t) \Vert_\mathcal{H} \label{eq: first maj norm X(t)} \\
& \phantom{\leq}\; + \left\Vert \int_0^t S(t-\tau) \left\{  \mathcal{A} B d(\tau) + U(\tau) \right\} \diff\tau \right\Vert_\mathcal{H} .  \nonumber
\end{align}
Let $m_r,M_R \in \mathbb{R}_+^*$ be the constants associated with the inequality (\ref{eq: Riesz basis - inequality}) for the Riesz basis formed by the eigenvectors $\phi_n$ of $\mathcal{A}_0$. Introducing $\kappa_0 = - \omega_0 > 0$ where $\omega_0$ is the growth bound of $S$, it is easy to see based on (\ref{eq: Riesz basis - inequality}) and (\ref{eq: Riesz-spectral operator - expression semigroup}) that
\begin{equation}\label{eq: first maj norm X(t) - term 1}
\Vert S(t) X_0 \Vert_\mathcal{H} \leq \sqrt{\dfrac{M_R}{m_R}} e^{-\kappa_0 t} \Vert X_0 \Vert_\mathcal{H} . 
\end{equation}
Similarly,
\begin{align}
& \left\Vert \int_0^t S(t-\tau) \left\{  \mathcal{A} B d(\tau) + U(\tau) \right\} \diff\tau \right\Vert_\mathcal{H} \nonumber \\
\leq & \sqrt{\dfrac{M_R}{m_R}} \int_0^t e^{-\kappa_0 (t-\tau)} \Vert \mathcal{A} B d(\tau) + U(\tau) \Vert_\mathcal{H} \diff\tau \nonumber\\
\leq & \dfrac{1}{\kappa_0} \sqrt{\dfrac{M_R}{m_R}} \left\{ \Vert \mathcal{A} B \Vert_{\mathcal{L}(\mathbb{K}^m,\mathcal{H})} \Vert d \Vert_{\mathcal{C}^0([0,t],\mathbb{K}^m)} + \Vert U \Vert_{\mathcal{C}^0([0,t];\mathcal{H})} \right\} \label{eq: first maj norm X(t) - term 2} .
\end{align}
It remains to evaluate $\Vert \alpha(t) \Vert_\mathcal{H}$. To do so, consider a basis $\mathcal{E} = (e_1,e_2,\ldots,e_m)$ of $\mathbb{K}^m$. We introduce $d_1,d_2,\ldots,d_m \in \mathcal{C}^2(\mathbb{R}_+;\mathbb{K})$ such that
\begin{equation*}
d = \sum\limits_{k=1}^{m} d_k e_k .
\end{equation*}
Based on this projection, one can get for all $\tau \in [0,t]$,
\begin{align*}
\vert \left< B d(\tau) , \psi_n \right>_\mathcal{H} \vert 
& = \left\vert \sum\limits_{k=1}^{m} d_k(\tau) \left< B e_k , \psi_n \right>_\mathcal{H} \right\vert \nonumber \\
& \leq \Vert d(\tau) \Vert_{\infty,\mathcal{E}} \sum\limits_{k=1}^{m} \vert \left< B e_k , \psi_n \right>_\mathcal{H} \vert \nonumber \\
& \leq c(\mathcal{E}) \Vert d(\tau) \Vert_{\mathbb{K}^m} \sum\limits_{k=1}^{m} \vert \left< B e_k , \psi_n \right>_\mathcal{H} \vert \nonumber  \\
& \leq c(\mathcal{E}) \Vert d \Vert_{C^0([0,t];\mathbb{K}^m)} \sum\limits_{k=1}^{m} \vert \left< B e_k , \psi_n \right>_\mathcal{H} \vert . \nonumber 
\end{align*}
Thus, for all $t \geq 0$, we have
\begin{align}
\vert \alpha_n(t) \vert 
\leq & \vert \lambda_n \vert \int_0^t e^{\operatorname{Re}\lambda_n (t-\tau)} \vert \left< B d(\tau) , \psi_n \right>_\mathcal{H}\vert \diff\tau \label{eq: estimation alpha_n} \\
\leq & \left\vert \dfrac{\lambda_n}{\operatorname{Re}\lambda_n} \right\vert c(\mathcal{E}) \Vert d \Vert_{C^0([0,t];\mathbb{K}^m)} \sum\limits_{k=1}^{m} \vert \left< B e_k , \psi_n \right>_\mathcal{H} \vert \nonumber \\
& \times \int_0^t - \operatorname{Re} \lambda_n \cdot e^{\operatorname{Re}\lambda_n (t-\tau)} \diff\tau \nonumber \\
\leq & \zeta \, c(\mathcal{E}) \Vert d \Vert_{C^0([0,t];\mathbb{K}^m)} \sum\limits_{k=1}^{m} \vert \left< B e_k , \psi_n \right>_\mathcal{H} \vert , \nonumber
\end{align}
where the eigenvalue constraints (\ref{eq: problem setting - confinement conditions}) have been used. We deduce that, for all $t \geq 0$,
\begin{equation*}
\vert \alpha_n(t) \vert^2 
\leq m \zeta^2 \, c(\mathcal{E})^2 \Vert d \Vert_{C^0([0,t];\mathbb{K}^m)}^2 \sum\limits_{k=1}^{m} \vert \left< B e_k , \psi_n \right>_\mathcal{H} \vert^2 .
\end{equation*}
From (\ref{eq: Riesz basis - inequality}), we deduce first that 
\begin{equation*}
\sum\limits_{n \in \mathbb{N}} \vert \alpha_n(t) \vert^2 
\leq \dfrac{m}{m_R} \zeta^2 \, c(\mathcal{E})^2 \Vert d \Vert_{C^0([0,t];\mathbb{K}^m)}^2 \sum\limits_{k=1}^{m} \Vert Be_k \Vert_\mathcal{H}^2 ,
\end{equation*}
and then, for all $t \geq 0$,
\begin{equation}\label{eq: first maj norm X(t) - term 3}
\Vert \alpha(t) \Vert_\mathcal{H} 
\leq \zeta \, c(\mathcal{E}) \sqrt{ m \dfrac{M_R}{m_R} \sum\limits_{k=1}^{m} \Vert Be_k \Vert_\mathcal{H}^2 } \Vert d \Vert_{C^0([0,t];\mathbb{K}^m)} .
\end{equation}
Substituting inequalities (\ref{eq: first maj norm X(t) - term 1}-\ref{eq: first maj norm X(t) - term 2}) and (\ref{eq: first maj norm X(t) - term 3}) into (\ref{eq: first maj norm X(t)}), we obtain the desired result (\ref{eq: th classical solutions - ISS}) with
\begin{equation*}
C_0 = \sqrt{\dfrac{M_R}{m_R}} ,
\qquad
C_2 = \dfrac{1}{\kappa_0} \sqrt{\dfrac{M_R}{m_R}} ,
\end{equation*}
\begin{equation}\label{eq: constant C_1}
C_1 = \sqrt{ \dfrac{M_R}{m_R} } \left\{ \dfrac{1}{\kappa_0} \Vert \mathcal{A} B \Vert_{\mathcal{L}(\mathbb{K}^m,\mathcal{H})} + \zeta \, c(\mathcal{E}) \sqrt{ m \sum\limits_{k=1}^{m} \Vert Be_k \Vert_\mathcal{H}^2 } \right\} .
\end{equation}
This concludes the proof.\qed

\begin{rem}\label{rem: novel spectral decomposition}
Equations (\ref{eq: projection of X into the Riesz basis}-\ref{eq: diff equation c_n}) actually hold true for the classical solutions associated with any Riesz-spectral boundary control system. Thus, the original Riesz-spectral boundary control system is equivalent to a countably infinite number of uncoupled ODEs describing the time domain evolution of the coefficients $c_n$ corresponding to the projection of the system trajectories onto the Riesz basis $\left\{ \phi_n , \; n \in \mathbb{N} \right\}$. The notable feature of the obtained ODEs (\ref{eq: diff equation c_n}) is that they provide a spectral decomposition that involves the boundary disturbance $d$ but avoids the occurrence of its time derivative $\dot{d}$. This important feature offers opportunities for feedback control design. We refer, e.g., to~\cite{coron2004global,coron2006global} for the global stabilization of heat and wave equations based on such a spectral decomposition but in the presence of a $\dot{d}$ term.
\end{rem}

\begin{rem}\label{rem: theorem 1 relax zeta}
In the proof of Theorem~\ref{th: ISS classical solutions}, the eigenvalue constraint $\zeta < \infty$ can be weakened to
\begin{equation}\label{eq: relaxed damping constraint}
\forall k \in \{1,\ldots,m\}, \;\;
\sum\limits_{n \geq 0} \left\vert \dfrac{\lambda_n}{\operatorname{Re}\lambda_n} \right\vert^2 \left\vert\left< B e_k , \psi_n \right>_\mathcal{H} \right\vert^2 < \infty .
\end{equation}
It is easy to see that the condition above does not depend on a specific selection of either the lifting operator $B$ (when $\omega_0 < 0$) or the basis $\mathcal{E} = (e_1,e_2,\ldots,e_m)$ of $\mathbb{K}^m$. In this case, the constant $C_1$ is given by
\begin{align}
C_1 & = \dfrac{1}{\kappa_0} \sqrt{ \dfrac{M_R}{m_R} } \Vert \mathcal{A} B \Vert_{\mathcal{L}(\mathbb{K}^m,\mathcal{H})} \label{eq: relaxed damping constraint - C1} \\
& \phantom{=}\, + c(\mathcal{E}) \sqrt{ m M_R \sum\limits_{k=1}^{m}  \sum\limits_{n \geq 0} \left\vert \dfrac{\lambda_n}{\operatorname{Re}\lambda_n} \right\vert^2 \left\vert\left< B e_k , \psi_n \right>_\mathcal{H} \right\vert^2 } . \nonumber
\end{align}
\end{rem}

\subsection{An energy-based interpretation for the constant related to the boundary perturbation in the ISS estimate}\label{subsec: energy-based ISS estimate}
The obtained expression (\ref{eq: constant C_1}) of the constant $C_1$ depends on the selected lifting operator $B$. However, the lifting operator provided by Definition~\ref{def: boundary control system} is not unique. The objective of this subsection is to provide a constructive definition of a constant $C_1$, independent of a specific selection of the lifting operator $B$, such that the ISS estimate (\ref{eq: th classical solutions - ISS}) holds true.

\begin{lem}\label{lemma: algebraic system}
Let $(\mathcal{A},\mathcal{B})$ be a boundary control system such that $0 \in \rho(\mathcal{A}_0)$. For any $e\in\mathbb{K}^m$, there exists a unique $X_e \in D(\mathcal{A})$ such that $\mathcal{A}X_e = 0$ and $\mathcal{B}X_e = e$. Furthermore, if $B$ is a lifting operator associated with $(\mathcal{A},\mathcal{B})$, then 
\begin{equation*}
X_e = Be -\mathcal{A}_0^{-1}\mathcal{A}Be .
\end{equation*}
\end{lem}

\textbf{Proof of Lemma~\ref{lemma: algebraic system}.} For the uniqueness part, by linearity, it is sufficient to check that $\mathcal{A}X = 0$ and $\mathcal{B}X = 0$ implies $X = 0$. But $\mathcal{B}X = 0$ implies $X \in D(\mathcal{A}_0)$ whence $\mathcal{A}_0 X = 0$. As $\mathcal{A}_0$ is injective, this yields $X=0$. For the existence part, let $B$ be a lifting operator associated with the boundary control system $(\mathcal{A},\mathcal{B})$ as provided by Definition~\ref{def: boundary control system}. Consider $\tilde{X}_e = X_e - B e$. Then for $X_e \in D(\mathcal{A})$,
\begin{align*}
& \;\; \mathcal{A}X_e = 0 \;\mathrm{and}\; \mathcal{B}X_e = e \\
\Leftrightarrow & \;\; \mathcal{A}\tilde{X}_e = -\mathcal{A}Be \;\mathrm{and}\; \mathcal{B}\tilde{X}_e = 0 \\
\Leftrightarrow & \;\; \mathcal{A}_0\tilde{X}_e = -\mathcal{A}Be \;\mathrm{and}\; \tilde{X}_e \in D(\mathcal{A}_0) \\
\Leftrightarrow & \;\; \tilde{X}_e = -\mathcal{A}_0^{-1}\mathcal{A}Be .
\end{align*} 
Thus, $X_e = Be -\mathcal{A}_0^{-1}\mathcal{A}Be$ is the unique solution. \qed

\begin{rem}\label{rem: xe sol stationnary problem}
The stationary trajectory $X(t) = X_e$ provided by Lemma~\ref{lemma: algebraic system} is the classical solution of the abstract boundary control system (\ref{def: boundary control system}) associated with the initial condition $X_0=X_e$, the constant boundary disturbance $d(t)=e$, and the zero distributed disturbance $U=0$.
\end{rem}

We can now introduce the main result of this section.

\begin{thm}\label{thm: constant C1 independant of B}
Let $\mathcal{E} = (e_1,e_2,\ldots,e_m)$ be a basis of $\mathbb{K}^m$. Under the assumptions of Theorem~\ref{th: ISS classical solutions}, the conclusion of the theorem holds true with constants $C_0,C_1,C_2>0$ involved in the ISS estimate (\ref{eq: th classical solutions - ISS}) given by
\begin{equation*}
C_0 = \sqrt{\dfrac{M_R}{m_R}} ,
\qquad
C_2 = \dfrac{1}{\kappa_0} \sqrt{\dfrac{M_R}{m_R}} ,
\end{equation*}
\begin{equation*}
C_1 = \zeta \, c(\mathcal{E}) \sqrt{ m \dfrac{M_R}{m_R} \sum\limits_{k=1}^{m} \Vert X_{e,k} \Vert_\mathcal{H}^2 } ,
\end{equation*}
where $X_{e,k} \in D(\mathcal{A})$ is, for all $k\in\{1,\ldots,m\}$, the unique solution of $\mathcal{A}X_{e,k} = 0$ and $\mathcal{B}X_{e,k} = e_k$.
\end{thm} 

\textbf{Proof of Theorem~\ref{thm: constant C1 independant of B}.} 
Consider the proof of Theorem~\ref{th: ISS classical solutions} up to Equation (\ref{eq: diff equation c_n}) included. For the basis $\mathcal{E} = (e_1,e_2,\ldots,e_m)$ of $\mathbb{K}^m$, let $d_1,d_2,\ldots,d_m \in \mathcal{C}^2(\mathbb{R}_+;\mathbb{K})$ be such that
\begin{equation*}
d = \sum\limits_{k=1}^{m} d_k e_k .
\end{equation*}
As $\mathcal{A}_0$ is a Riesz-spectral operator with $\omega_0 < 0$, we have $0 \in \rho(\mathcal{A}_0)$. Based on Lemma~\ref{lemma: algebraic system}, let, for any $k\in\{1,\ldots,m\}$, $X_{e,k} \in D(\mathcal{A})$ be the unique solution of $\mathcal{A}X_{e,k} = 0$ and $\mathcal{B}X_{e,k} = e_k$. In particular, an explicit expression is given by $X_{e,k} = Be_k -\mathcal{A}_0^{-1}\mathcal{A}Be_k$. Introducing $\tilde{X}_{e,k} \triangleq X_{e,k} - Be_k = -\mathcal{A}_0^{-1}\mathcal{A}Be_k$, this yields $\tilde{X}_{e,k} \in D(\mathcal{A}_0)$. Thus we have,
\begin{align*}
\lambda_n \left< X_{e,k} , \psi_n \right>_\mathcal{H}
& = \left< X_{e,k} , \overline{\lambda_n} \psi_n \right>_\mathcal{H} \\
& = \left< X_{e,k} , \mathcal{A}_0^* \psi_n \right>_\mathcal{H} \\
& = \left< \tilde{X}_{e,k} , \mathcal{A}_0^* \psi_n \right>_\mathcal{H} + \left< B e_k , \mathcal{A}_0^* \psi_n \right>_\mathcal{H} \\
& = \left< \mathcal{A}_0 \tilde{X}_{e,k} , \psi_n \right>_\mathcal{H} + \lambda_n \left< B e_k , \psi_n \right>_\mathcal{H} \\
& = - \left< \mathcal{A}Be_k , \psi_n \right>_\mathcal{H} + \lambda_n \left< B e_k , \psi_n \right>_\mathcal{H} .
\end{align*}
Introducing that
\begin{equation}\label{eq: project boundary disturbance over stationnary sol}
D = \sum\limits_{k=1}^{m} d_k X_{e,k} \in C^0(\mathbb{R}_+;D(\mathcal{A})) \cap C^2(\mathbb{R}_+;\mathcal{H}) , 
\end{equation}
it follows that
\begin{align*}
& - \lambda_n \left< B d(t) , \psi_n \right>_\mathcal{H} + \left< \mathcal{A} B d(t) , \psi_n \right>_\mathcal{H} \\
& = \sum\limits_{k=1}^{m} d_k(t) \left\{ - \lambda_n \left< B e_k , \psi_n \right>_\mathcal{H} + \left< \mathcal{A} B e_k , \psi_n \right>_\mathcal{H} \right\} \\
& = - \lambda_n \sum\limits_{k=1}^{m} d_k(t) \left< X_{e,k} , \psi_n \right>_\mathcal{H} \\
& = - \lambda_n \left< D(t) , \psi_n \right>_\mathcal{H} .
\end{align*}
Thus, we get from (\ref{eq: diff equation c_n}) that for all $t \geq 0$,
\begin{equation*}
\dot{c}_n(t) 
= \lambda_n c_n(t) - \lambda_n \left< D(t) , \psi_n \right>_\mathcal{H} \nonumber
+ \left< U(t) , \psi_n \right>_\mathcal{H} ,
\end{equation*}
and that for all $t \geq 0$,
\begin{align}
c_n(t) 
& = e^{\lambda_n t} c_n(0) - \lambda_n \int_0^t e^{\lambda_n (t-\tau)} \left< D(\tau) , \psi_n \right>_\mathcal{H} \diff\tau \nonumber \\
& \phantom{=}\; + \int_0^t e^{\lambda_n (t-\tau)} \left< U(\tau) , \psi_n \right>_\mathcal{H} \diff\tau . \nonumber
\end{align}
The first and third terms on the right hand side of the equation above have been estimated in the proof of Theorem~\ref{th: ISS classical solutions}. This procedure yields the constants $C_0$ and $C_2$. The second term can be treated with the same procedure that one employed for (\ref{eq: coeff alpha(t)}) via the projection (\ref{eq: project boundary disturbance over stationnary sol}). This also provides the estimate of $C_1$. \qed

\begin{rem}
The constant $C_1$ given by Theorem~\ref{thm: constant C1 independant of B} depends on the energy $\Vert X_{e,k} \Vert_\mathcal{H}$ of $m$ linearly independent\footnote{Directly follows from the definition of $X_{e,k}$ and the fact that $\mathcal{E}$ is a basis.} stationary solutions $X_{e,k}$ of the abstract boundary control system that are associated with the constant boundary perturbations $d(t) = e_k$ and the zero distributed disturbance $U=0$. 
\end{rem}

\begin{rem}\label{rem: fading memory estimates}
One of the main applications of ISS estimates relies in the establishment of small gain conditions for assessing the stability of interconnected systems~\cite[Sec.~8.1]{karafyllis2019preview}. This requires the conversion of the ISS estimates into fading memory estimates, e.g., by means of a conversion lemma~\cite[Chap.~7]{karafyllis2019preview}. Nevertheless, in the context of Theorems~\ref{th: ISS classical solutions} and~\ref{thm: constant C1 independant of B}, it is actually not necessary to resort to such a conversion lemma because a slight adaptation of the associated proofs (specifically, estimates (\ref{eq: first maj norm X(t) - term 2}-\ref{eq: estimation alpha_n})) directly shows that the following fading memory estimate holds true for all $\epsilon \in [0,1)$ and all $t \geq 0$, 
\begin{align*}
\left\Vert X(t) \right\Vert_\mathcal{H}
\leq &
C_0 e^{-\kappa_0 t} \left\Vert X_0 \right\Vert_\mathcal{H} 
+ \dfrac{C_1}{1-\epsilon} \Vert e^{- \epsilon \kappa_0 (t - (\cdot))} d \Vert_{\mathcal{C}^0([0,t],\mathbb{K}^m)} \nonumber \\
& + \dfrac{C_2}{1-\epsilon} \Vert e^{- \epsilon \kappa_0 (t - (\cdot))} U \Vert_{\mathcal{C}^0([0,t];\mathcal{H})} .
\end{align*}
\end{rem}

\section{Concept of Weak Solutions and ISS Property}\label{sec: ISS weak solutions}

Assume that a given boundary control system\footnote{Not necessarily a Riesz-spectral one.} $(\mathcal{A},\mathcal{B})$ satisfies an ISS estimate\footnote{In the sense provided later in Theorem~\ref{th: extended main result}.} with respect to classical solutions. The objective of this section is to extend such an ISS estimate to every initial condition $X_0 \in \mathcal{H}$ and every disturbance $d \in \mathcal{C}^0(\mathbb{R}_+;\mathbb{K}^m)$ and $U \in \mathcal{C}^0(\mathbb{R_+};\mathcal{H})$. To do so, we introduce a concept of weak solutions that extends to abstract boundary control systems the concept of weak solutions originally introduced for infinite dimensional nonhomogeneous Cauchy problems in~\cite{ball1977strongly} and further investigated in~\cite[Def.~3.1.6, A.5.29]{Curtain2012}. This represents an alternative to the traditional concept of mild solutions defined within the extrapolation space~\cite{emirsjlow2000pdes}.

\subsection{Definition of weak solutions}
To motivate the definition of weak solutions, let us consider first $X$ a classical solution associated with $(X_0,d,U)$. Let $T>0$ and a function $z \in \mathcal{C}^0([0,T];D(\mathcal{A}_0^*))$ be arbitrarily given. As $X(t)-Bd(t) \in D(\mathcal{A}_0)$, we obtain for all $t \in [0,T]$,
\begin{align*}
\left< \dfrac{\mathrm{d} X}{\mathrm{d} t}(t),z(t)\right>_\mathcal{H} 
& = \left< \mathcal{A}X(t)+U(t),z(t)\right>_\mathcal{H} \\
& = \left< \mathcal{A}_0 \left\{ X(t) - B d(t) \right\} , z(t) \right>_\mathcal{H} \\
& \phantom{=}\; + \left< \mathcal{A} B d(t) , z(t) \right>_\mathcal{H} + \left< U(t) , z(t) \right>_\mathcal{H} \\
& = \left< X(t) , \mathcal{A}_0^* z(t) \right>_\mathcal{H} - \left< B d(t) , \mathcal{A}_0^* z(t) \right>_\mathcal{H} \\
& \phantom{=}\; + \left< \mathcal{A} B d(t) , z(t) \right>_\mathcal{H} + \left< U(t) , z(t) \right>_\mathcal{H} .
\end{align*}

Assuming now that $z \in \mathcal{C}^0([0,T];D(\mathcal{A}_0^*)) \cap \mathcal{C}^1([0,T];\mathcal{H})$, integration by parts gives
\begin{align*}
\int_0^T \left< \dfrac{\mathrm{d} X}{\mathrm{d} t}(t) , z(t) \right>_\mathcal{H} \diff t 
& =
\left< X(T) , z(T) \right>_\mathcal{H}
- \left< X_0 , z(0) \right>_\mathcal{H} \\
& \phantom{=}\;  - \int_0^T \left< X(t) , \dfrac{\mathrm{d} z}{\mathrm{d} t}(t) \right>_\mathcal{H} \diff t .
\end{align*}

Based on the two identities above, we introduce the following definition.

\begin{defn}[Weak solutions]\label{def: weak solution}
Let $(\mathcal{A},\mathcal{B})$ be a boundary control system. For $X_0 \in \mathcal{H}$ and disturbances $d\in\mathcal{C}^0(\mathbb{R}_+;\mathbb{K}^m)$ and $U\in\mathcal{C}^0(\mathbb{R}_+;\mathcal{H})$, we say that $X \in \mathcal{C}^0(\mathbb{R}_+;\mathcal{H})$ is a weak solution of the abstract boundary control system (\ref{def: boundary control system}) associated with $(X_0,d,U)$ if for all $T>0$ and for all $z \in \mathcal{C}^0([0,T];D(\mathcal{A}^*_0)) \cap \mathcal{C}^1([0,T];\mathcal{H})$ such that $\mathcal{A}_0^* z \in \mathcal{C}^0([0,T];\mathcal{H})$ and $z(T)=0$ (such a function $z$ is called a test function over $[0,T]$), the following equality holds true: 
\begin{align}
& \int_0^T \left< X(t), \mathcal{A}_0^* z(t) + \dfrac{\mathrm{d} z}{\mathrm{d} t}(t) \right>_\mathcal{H} \diff t \nonumber \\
& = - \left< X_0 , z(0) \right>_\mathcal{H} 
+ \int_0^T \left< B d(t) , \mathcal{A}_0^* z(t) \right>_\mathcal{H} \diff t \label{eq: def weak solution} \\
& \phantom{=}\; - \int_0^T \left< \mathcal{A} B d(t) , z(t) \right>_\mathcal{H} \diff t - \int_0^T \left< U(t) , z(t) \right>_\mathcal{H} \diff t \nonumber ,
\end{align}
where $B$ is an arbitrary lifting operator associated with $(\mathcal{A},\mathcal{B})$.
\end{defn}

\begin{rem}
Definition~\ref{def: weak solution} is relevant because for any $T > 0$, the set of test functions over $[0,T]$ is not reduced to only the zero function. Indeed, the function $z$ defined for all $t \in [0,T]$ by $z(t) = (t-T) z_0$ with $z_0 \in D(\mathcal{A}_0^*) \backslash\{0\}$ is obviously a non zero test function over $[0,T]$.
\end{rem}

\begin{rem}\label{rmk: classical implies weak}
The definition of a weak solution for the abstract system (\ref{def: boundary control system}) is compatible with the notion of classical solution. Indeed, the developments preliminary to the introduction of Definition~\ref{def: weak solution} show that a classical solution is also a weak solution.
\end{rem}

\begin{rem}
At first sight, the right hand side of (\ref{eq: def weak solution}) depends on the selected lifting operator $B$, making it necessary to specify the selected lifting operator $B$ when saying that a trajectory $X$ is a weak solution associated with $(X_0,d,U)$. However, Definition~\ref{def: weak solution} implicitly claims that a weak solution is actually independent of the selected lifting operator $B$. This directly follows from the fact that if $B$ and $\tilde{B}$ are two lifting operators associated with $(\mathcal{A},\mathcal{B})$, then $\hat{B} \triangleq B - \tilde{B}$ satisfies $\mathrm{R}(\hat{B}) \subset D(\mathcal{A})$ and $\mathcal{B}\hat{B} = \mathcal{B}B - \mathcal{B}\tilde{B} = I_{\mathbb{K}^m} - I_{\mathbb{K}^m} = 0$, i.e., $\mathrm{R}(\hat{B}) \subset D(\mathcal{A}_0)$. Thus we obtain, $\left< \hat{B} d(t) , \mathcal{A}_0^* z(t) \right>_\mathcal{H} = \left< \mathcal{A}_0 \hat{B} d(t) , z(t) \right>_\mathcal{H} = \left< \mathcal{A} \hat{B} d(t) , z(t) \right>_\mathcal{H}$, from which we deduce that
\begin{align*}
& \left< B d(t) , \mathcal{A}_0^* z(t) \right>_\mathcal{H} 
- \left< \mathcal{A} B d(t) , z(t) \right>_\mathcal{H} \\
& \phantom{======} = \left< \tilde{B} d(t) , \mathcal{A}_0^* z(t) \right>_\mathcal{H}
- \left< \mathcal{A} \tilde{B} d(t) , z(t) \right>_\mathcal{H} .
\end{align*}
This equality shows that the right hand side of (\ref{eq: def weak solution}) remains unchanged when switching between different lifting operators $B$ associated with $(\mathcal{A},\mathcal{B})$. So, it indeed makes sense to discuss about weak solutions without mentioning a particular lifting operator $B$ associated with $(\mathcal{A},\mathcal{B})$.
\end{rem} 

\begin{rem}
Note that the concept of weak solution does not require that the initial condition satisfies the boundary condition $\mathcal{B}X_0 = d(0)$. Such an algebraic condition is not even well defined when $X_0 \notin D(\mathcal{B})$.
\end{rem}

\subsection{Properties of weak solutions}

When defining a notion of a weak solution, it is generally desirable to preserve the uniqueness of the solution.

\begin{lem}\label{lemma: weak solution uniqueness}
Let $(\mathcal{A},\mathcal{B})$ be a boundary control system such that the disturbance-free operator $\mathcal{A}_0$ is injective. Then, for any given $X_0 \in \mathcal{H}$ and disturbances $d\in\mathcal{C}^0(\mathbb{R}_+;\mathbb{K}^m)$ and $u\in\mathcal{C}^0(\mathbb{R}_+;\mathcal{H})$, there exists at most one weak solution associated with $(X_0,d,U)$ of the abstract system (\ref{def: boundary control system}).
\end{lem}

The proof of Lemma~\ref{lemma: weak solution uniqueness} is in Annex~\ref{annex: proof lemma uniqueness weak solutions}. We deduce from Remark~\ref{rmk: classical implies weak} and Lemma~\ref {lemma: weak solution uniqueness} the following result.

\begin{cor}\label{cor: equiv classical and weak solutions for regular IC and disturbances}
Let $(\mathcal{A},\mathcal{B})$ be a boundary control system such that the disturbance-free operator $\mathcal{A}_0$ is injective. For any initial condition $X_0 \in D(\mathcal{A})$ and disturbances $d \in \mathcal{C}^2(\mathbb{R}_+;\mathbb{K}^m)$ and $U \in \mathcal{C}^1(\mathbb{R}_+;\mathcal{H})$ such that $\mathcal{B}X_0 = d(0)$, the concepts of classical and weak solutions coincide. More specifically, the two following statements are equivalent:
\begin{enumerate}
\item $X$ is a classical solution associated with $(X_0,d,U)$;
\item $X$ is a weak solution associated with $(X_0,d,U)$.
\end{enumerate}
\end{cor}

When $X$ is a classical solution of the abstract boundary control system (\ref{def: boundary control system}) associated with $(X_0,d,U)$, we have by Definition~\ref{def: classical solution} that $X(0)=X_0$. Such an initial condition is not explicitly imposed in the Definition~\ref{def: weak solution} of a weak solution. However, it is a consequence of (\ref{eq: def weak solution}) as shown by the following lemma.

\begin{lem}\label{lem: weak solution initial condition}
Let $(\mathcal{A},\mathcal{B})$ be a boundary control system. Under the terms of Definition~\ref{def: weak solution}, assume that $X$ is a weak solution associated with $(X_0,d,U)$. Then $X(0)=X_0$.
\end{lem}

\textbf{Proof of Lemma~\ref{lem: weak solution initial condition}.} Let $z \in \mathcal{C}^0(\mathbb{R}_+;D(\mathcal{A}^*_0)) \cap \mathcal{C}^1(\mathbb{R}_+;\mathcal{H})$ such that $\mathcal{A}_0^* z \in \mathcal{C}^0(\mathbb{R}_+;\mathcal{H})$. Then for all $T > 0$, $\tilde{z}_T \triangleq z - z(T)$ is a test function over $[0,T]$. As $\mathrm{d} z(T) / \mathrm{d}t = 0$, (\ref{eq: def weak solution}) gives for all $T > 0$,
\begin{align}
& \dfrac{1}{T} \int_0^T \left< X(t), \mathcal{A}_0^* z(t) - \mathcal{A}_0^* z(T) \right>_\mathcal{H} \diff t \nonumber \\
& + \dfrac{1}{T} \int_0^T \left< X(t), \dfrac{\mathrm{d} z}{\mathrm{d} t}(t) \right>_\mathcal{H} \diff t \nonumber \\
& = \left< X_0 , \dfrac{z(T) - z(0)}{T} \right>_\mathcal{H} \label{eq: X(0)=X_0 - finite difference} \\
& \phantom{=}\; + \dfrac{1}{T} \int_0^T \left< B d(t) , \{ \mathcal{A}_0^* z(t) - \mathcal{A}_0^* z(T) \} \right>_\mathcal{H} \diff t \nonumber \\
& \phantom{=}\; - \dfrac{1}{T} \int_0^T \left< \mathcal{A} B d(t) , z(t) - z(T) \right>_\mathcal{H} \diff t \nonumber \\
& \phantom{=}\; - \dfrac{1}{T} \int_0^T \left< U(t) , z(t) - z(T) \right>_\mathcal{H} \diff t . \nonumber 
\end{align}
It is straightforward to show that for any $f,g \in \mathcal{C}^0(\mathbb{R}_+;\mathcal{H})$,
\begin{equation*}
\underset{T \rightarrow 0^+}{\lim} \; \dfrac{1}{T} \int_0^T \left< f(t), g(t) - g(T) \right>_\mathcal{H} \diff t = 0.
\end{equation*}
Thus, based on the regularity assumptions given in Definition~\ref{def: weak solution}, we deduce, by letting $T \rightarrow 0^+$ in (\ref{eq: X(0)=X_0 - finite difference}),  
\begin{equation*}
\left< X(0), \dfrac{\mathrm{d} z}{\mathrm{d} t}(0) \right>_\mathcal{H}
= \left< X_0 , \dfrac{\mathrm{d} z}{\mathrm{d} t}(0) \right>_\mathcal{H} .
\end{equation*} 
Taking in particular $z(t) = t z_0$ with $z_0 \in D(\mathcal{A}_0^*)$, we get that $\left< X(0) - X_0, z_0 \right>_\mathcal{H} = 0$ holds true for all $z_0 \in D(\mathcal{A}_0^*)$. As $\overline{D(\mathcal{A}_0^*)} = \mathcal{H}$, we deduce that $X(0)=X_0$. \qed

Thus, for a weak solution $X$ associated with $(X_0,d,u)$, it makes sense to say that $X_0$ is the initial condition of the system trajectory.

\subsection{Existence of weak solutions and extension of ISS estimates}

The proof of the following theorem is in Annex~\ref{annex: proof thm extended main result}.

\begin{thm}\label{th: extended main result}
Let $(\mathcal{A},\mathcal{B})$ be a boundary control system. Assume that there exist $\beta \in \mathcal{KL}$ and $\gamma_1,\gamma_2\in\mathcal{K}$ such that for any initial condition $X_0 \in D(\mathcal{A})$ and any disturbances $d \in \mathcal{C}^2(\mathbb{R}_+;\mathbb{K}^m)$ and $U \in \mathcal{C}^1(\mathbb{R}_+;\mathcal{H})$ such that $\mathcal{B}X_0 = d(0)$, the classical solution $X$ of the abstract boundary control system (\ref{def: boundary control system}) associated with $(X_0,d,U)$ satisfies for all $t \geq 0$,
\begin{align}
\left\Vert X(t) \right\Vert_\mathcal{H}
\leq &
\beta(\left\Vert X_0 \right\Vert_\mathcal{H},t) 
+ \gamma_1 \left(\Vert d \Vert_{\mathcal{C}^0([0,t],\mathbb{K}^m)} \right) \nonumber \\
& + \gamma_2 \left( \Vert U \Vert_{\mathcal{C}^0([0,t];\mathcal{H})} \right) . \label{eq: th weak solutions - ISS}
\end{align}
Then, for any initial condition $X_0 \in \mathcal{H}$, and any disturbances $d\in\mathcal{C}^0(\mathbb{R}_+;\mathbb{K}^m)$ and $U\in\mathcal{C}^0(\mathbb{R}_+;\mathcal{H})$:
\begin{enumerate}
\item the abstract boundary control system (\ref{def: boundary control system}) has a unique weak solution $X \in \mathcal{C}^0(\mathbb{R}_+;\mathcal{H})$ associated with $(X_0,d,U)$; 
\item this weak solution satisfies the ISS estimate (\ref{eq: th weak solutions - ISS}) for all $t \geq 0$.
\end{enumerate} 
\end{thm} 

\begin{rem}
Theorem~\ref{th: extended main result} ensures that the ISS property holds for weak solutions if and only if it holds for classical solutions. Thus, in the study of a given abstract boundary control system, it is actually sufficient to study the ISS property for classical solutions associated with sufficiently smooth disturbances (from the density argument used in the proof of Theorem~\ref{th: extended main result}, the study can even be restricted to the only disturbances of class $C^\infty$) to conclude that the ISS property holds for weak solutions.
\end{rem}

We directly deduce from Theorem~\ref{th: extended main result} the following extension of Theorems~\ref{th: ISS classical solutions} and~\ref{thm: constant C1 independant of B} for Riesz-spectral operators.

\begin{cor}
Let $(\mathcal{A},\mathcal{B})$ be a Riesz-spectral boundary control system such that the eigenvalue constraints (\ref{eq: problem setting - confinement conditions}) hold true\footnote{The constraint $\zeta < \infty$ can be relaxed to (\ref{eq: relaxed damping constraint}). In that case, the constant $C_1$ is given by (\ref{eq: relaxed damping constraint - C1}).}. For every initial condition $X_0 \in \mathcal{H}$, and every disturbance $d \in \mathcal{C}^0(\mathbb{R}_+;\mathbb{K}^m)$ and $U \in \mathcal{C}^0(\mathbb{R}_+;\mathcal{H})$, the abstract system (\ref{def: boundary control system}) has a unique weak solution $X \in \mathcal{C}^0(\mathbb{R}_+;\mathcal{H})$ associated with $(X_0,d,U)$. Furthermore, $X$ satisfies the ISS estimate (\ref{eq: th classical solutions - ISS}) with constants $\kappa_0,C_0,C_1,C_2$ given by either Theorem~\ref{th: ISS classical solutions} or~\ref{thm: constant C1 independant of B}.  
\end{cor}

Based on the existence and uniqueness of weak solutions, and by linearity of (\ref{eq: def weak solution}), we can state the following result. 

\begin{cor}\label{cor: weak solution - linearity}
Assume that the assumptions of Theorem~\ref{th: extended main result} hold true. For $i\in\{1,2\}$, let $X_i$ be the weak solution associated with $(X_{i,0},d_i,U_i)$. Then, for all $\alpha,\beta\in\mathbb{K}$, $\alpha X_1 + \beta X_2$ is the unique weak solution associated with $(\alpha X_{1,0} + \beta X_{2,0},\alpha d_{1} + \beta d_{2},\alpha U_{1} + \beta U_{2})$. 
\end{cor}

\subsection{Mild solutions and semigroup property}

\subsubsection{Compatibility with the concept of mild solutions}
When the boundary disturbance satisfies the additional regularity assumption $d \in \mathcal{C}^1(\mathbb{R}_+;\mathbb{K}^m)$, a classical approach for extending the concept of classical solutions is to define $X$ provided by
\begin{align}
X(t) & =  S(t) (X_0 - Bd(0)) + Bd(t) \label{eq: explicit expression classical solution} \\
& \phantom{=}\; + \int_0^t S(t-\tau) \left\{ - B\dot{d}(\tau) +  \mathcal{A}Bd(\tau) + U(\tau) \right\} \diff\tau \nonumber
\end{align}
as the mild solution associated with $(X_0,d,U)$. In particular, any classical solution $X$ satisfies the identity (\ref{eq: explicit expression classical solution}). The following result shows that such an approach is compatible with the concept of weak solutions.

\begin{thm}\label{thm: explicit form weak solution d of class C1}
Let $(\mathcal{A},\mathcal{B})$ be a boundary control system such that the assumptions of Theorem~\ref{th: extended main result} hold true. Under the terms of Definition~\ref{def: weak solution}, let $X$ be the weak solution associated with $(X_0,d,U)$. Assume that the boundary disturbance satisfies the extra regularity assumption $d \in \mathcal{C}^1(\mathbb{R}_+;\mathbb{K}^m)$. Then $X$ is also a mild solution in the sense that (\ref{eq: explicit expression classical solution}) holds true for all $t \geq 0$.
\end{thm}

\textbf{Proof of Theorem~\ref{thm: explicit form weak solution d of class C1}.}
Let $X_0\in\mathcal{H}$, $d \in \mathcal{C}^1(\mathbb{R}_+;\mathbb{K}^m)$, $U \in \mathcal{C}^0(\mathbb{R}_+;\mathcal{H})$, and $T > 0$ be arbitrarily given. We denote by $X$ the weak solution associated with $(X_0,d,U)$. As $\mathcal{C}^2([0,T];\mathbb{K}^m)$ is dense in $\mathcal{C}^1([0,T];\mathbb{K}^m)$ endowed with the usual norm $\Vert f \Vert_{\mathcal{C}^1([0,T];\mathbb{K}^m)} = \Vert f \Vert_{\mathcal{C}^0([0,T];\mathbb{K}^m)} + \Vert \dot{f} \Vert_{\mathcal{C}^0([0,T];\mathbb{K}^m)}$, we can select, similarly to the proof of Theorem~\ref{th: extended main result}, approximating sequences $(X_{0,n})_n \in D(\mathcal{A})^{\mathbb{N}}$, $(d_{n})_n \in \mathcal{C}^2([0,T];\mathbb{K}^m)^\mathbb{N}$, and $(U_{n})_n \in \mathcal{C}^1([0,T];\mathcal{H})^\mathbb{N}$ such that $\mathcal{B}X_{0,n} = d_n(0)$, $X_{0,n} \underset{n \rightarrow +\infty}{\longrightarrow} X_0$, 
\begin{equation*}
\Vert d_n - d \Vert_{\mathcal{C}^0([0,T];\mathbb{K}^m)} \underset{n \rightarrow +\infty}{\longrightarrow} 0 ,
\end{equation*}
\begin{equation*}
\Vert \dot{d}_n - \dot{d} \Vert_{\mathcal{C}^0([0,T];\mathbb{K}^m)} \underset{n \rightarrow +\infty}{\longrightarrow} 0 ,
\end{equation*}
\begin{equation*}
\Vert U_n - U \Vert_{\mathcal{C}^0([0,T];\mathcal{H})} \underset{n \rightarrow +\infty}{\longrightarrow} 0 .
\end{equation*}
We denote by $X_n$ the unique classical solution of the abstract system (\ref{def: boundary control system}) over $[0,T]$ associated with $(X_{0,n},d_n,U_n)$. From (\ref{eq: explicit expression classical solution}), we obtain for all $t \in [0,T]$ and all $n \in \mathbb{N}$,
\begin{align}
X_n(t) & =  S(t) (X_{0,n} - Bd_n(0)) + Bd_n(t) \label{eq: explicit expression classical solution - proof mild sol} \\
& \phantom{=}\; + \int_0^t S(t-\tau) \left\{ - B\dot{d}_n(\tau) +  \mathcal{A}Bd_n(\tau) + U_n(\tau) \right\} \diff\tau . \nonumber 
\end{align}
Furthermore, from the proof of Theorem~\ref{th: extended main result}, we know that  $\Vert X_n - X \Vert_{\mathcal{C}^0([0,T];\mathcal{H})} \underset{n \rightarrow +\infty}{\longrightarrow} 0$. Thus, by letting $n \rightarrow +\infty$ in (\ref{eq: explicit expression classical solution - proof mild sol}), we obtain that (\ref{eq: explicit expression classical solution}) holds true for all $t \in [0,T]$. As $T > 0$ has been arbitrarily chosen, this concludes the proof.\qed

\subsubsection{Semigroup property}
It is well known that the classical/mild solutions of the abstract system (\ref{def: boundary control system}) satisfy the semigroup property in the sense that if $X$ is the classical/mild solution associated with $(X_0,d,U)$, then $X(\cdot+t_0)$ is the classical/mild solution associated with  $(X(t_0),d(\cdot+t_0),U(\cdot+t_0))$ for any $t_0 > 0$. The following result shows that this semigroup property extends to the concept of weak solutions.

\begin{thm}\label{thm: time-invariant}
Let $(\mathcal{A},\mathcal{B})$ be a boundary control system such that the assumptions of Theorem~\ref{th: extended main result} hold true. Let $X$ be the weak solution associated with an initial condition $X_0 \in \mathcal{H}$, a boundary disturbance $d\in\mathcal{C}^0(\mathbb{R}_+;\mathbb{K}^m)$, and a distributed disturbance $U\in\mathcal{C}^0(\mathbb{R}_+;\mathcal{H})$. Then, for any $t_0 > 0$, $X(\cdot+t_0)$ is the weak solution associated with  $(X(t_0),d(\cdot+t_0),U(\cdot+t_0))$.
\end{thm}

\textbf{Proof of Theorem~\ref{thm: time-invariant}}
Note first that, based on the developments preliminary to the introduction of Definition~\ref{def: weak solution}, we have for any classical solution $X$ associated with $(X_0,d,U)$, any $T > 0$, and any $z \in \mathcal{C}^0([0,T];D(\mathcal{A}^*_0)) \cap \mathcal{C}^1([0,T];\mathcal{H})$ such that\footnote{We do not impose here the condition $z(T)=0$ as in the case of the test functions.} $\mathcal{A}_0^* z \in \mathcal{C}^0([0,T];\mathcal{H})$,
\begin{align}
& \int_0^T \left< X(t), \mathcal{A}_0^* z(t) + \dfrac{\mathrm{d} z}{\mathrm{d} t}(t) \right>_\mathcal{H} \diff t \nonumber \\
& = \left< X(T) , z(T) \right>_\mathcal{H} - \left< X_0 , z(0) \right>_\mathcal{H} \label{eq: weak solution - non zero z(T)} \\
& \phantom{=}\; + \int_0^T \left< B d(t) , \mathcal{A}_0^* z(t) \right>_\mathcal{H} \diff t \nonumber \\
& \phantom{=}\; - \int_0^T \left< \mathcal{A} B d(t) , z(t) \right>_\mathcal{H} \diff t - \int_0^T \left< U(t) , z(t) \right>_\mathcal{H} \diff t \nonumber .
\end{align}
By resorting to the same density argument as the one employed in Step~4 of the proof of Theorem~\ref{th: extended main result}, this yields that any weak solution $X$ associated with $(X_0,d,U)$ also satisfies (\ref{eq: weak solution - non zero z(T)}) for all $T > 0$ and all $z \in \mathcal{C}^0([0,T];D(\mathcal{A}^*_0)) \cap \mathcal{C}^1([0,T];\mathcal{H})$ such that $\mathcal{A}_0^* z \in \mathcal{C}^0([0,T];\mathcal{H})$.

Now, let $X$ be the weak solution associated with $(X_0,d,U)$. Let $t_0,T > 0$ and a test function $\hat{z} \in \mathcal{C}^0([0,T];D(\mathcal{A}^*_0)) \cap \mathcal{C}^1([0,T];\mathcal{H})$ over $[0,T]$ be arbitrarily given. We define the test function $\tilde{z} \in \mathcal{C}^0([0,t_0+T];D(\mathcal{A}^*_0)) \cap \mathcal{C}^1([0,t_0+T];\mathcal{H})$ over $[0,t_0+T]$ as $\tilde{z} = \left.\varphi(\cdot - t_0)\right\vert_{[0,t_0+T]}$ where $\varphi \in \mathcal{C}^0(\mathbb{R};D(\mathcal{A}^*_0)) \cap \mathcal{C}^1(\mathbb{R};\mathcal{H})$ is given for all $t\in[0,T]$ and all $k \in \mathbb{Z}$ by $\varphi(t+2kT)=\hat{z}(t)-2k\hat{z}(0)$ and $\varphi(t+(2k+1)T)=-\hat{z}(T-t)-2k\hat{z}(0)$. In particular we have $\tilde{z}(t_0+T)=\varphi(T)=\hat{z}(T)=0$ and for all $t \in [0,T]$, $\tilde{z}(t+t_0)=\varphi(t)=\hat{z}(t)$. By using (\ref{eq: weak solution - non zero z(T)}) once with $z = \left.\tilde{z}\right\vert_{[0,t_0]}$ over $[0,t_0]$ and once with $z = \tilde{z}$ over $[0,t_0+T]$, we obtain after a change of variable:
\begin{align*}
& \int_{0}^{T} \left< X(t+t_0), \mathcal{A}_0^* \hat{z}(t) + \dfrac{\mathrm{d} \hat{z}}{\mathrm{d} t}(t) \right>_\mathcal{H} \diff t \nonumber \\
& = - \left< X(t_0) , \hat{z}(0) \right>_\mathcal{H} 
+ \int_{0}^{T} \left< B d(t+t_0) , \mathcal{A}_0^* \hat{z}(t) \right>_\mathcal{H} \diff t \nonumber \\
& \phantom{=}\; - \int_{0}^{T} \left< \mathcal{A} B d(t+t_0) , \hat{z}(t) \right>_\mathcal{H} \diff t \nonumber \\
& \phantom{=}\; - \int_{0}^{T} \left< U(t+t_0) , \hat{z}(t) \right>_\mathcal{H} \diff t \nonumber .
\end{align*}
As $T$, $\hat{z}$ and $t_0$ have been arbitrarily selected, it follows from Definition~\ref{def: weak solution} that for all $t_0>0$, $X(\cdot+t_0)$ is the weak solution associated with $(X(t_0),d(\cdot+t_0),U(\cdot+t_0))$.
\qed

%

\section{Application}\label{sec: application}
Among the examples of applications, one can find 1D parabolic PDEs~\cite{karafyllis2016iss} and a flexible damped string~\cite{lhachemi2018input}. In this section, we detail another example of structural vibrations, namely a damped Euler-Bernoulli beam.

We denote by $L^2(0,1)$ and $H^m(0,1)$ the set of square (Lebesgue) integrable functions over $(0,1)$ and the usual Sobolev space of order $m$ over $(0,1)$, respectively. We also introduce $H_0^1(0,1) \triangleq \{ f \in H^1(0,1) \; : \; f(0) = f(1) = 0 \}$.

For a function $f : \mathbb{R}_+ \rightarrow L^2(0,1)$, we denote, with a slight abuse of notation, $f(t,\xi) \triangleq \left[f(t)\right](\xi)$. When $f\in\mathcal{C}^1(\mathbb{R}_+;\mathcal{H})$, we denote $\dfrac{\mathrm{d} f}{\mathrm{d} t}(t,\xi) \triangleq \left[ \dfrac{\mathrm{d} f}{\mathrm{d} t}(t) \right](\xi)$. Finally, when $f : \mathbb{R}_+ \rightarrow H^1(0,1)$, we denote $f'(t,\xi) \triangleq \left[ f(t) \right]'(\xi)$.

\subsection{Damped Euler-Bernoulli beam}\label{subsec: damped beam}
We consider a damped Euler-Bernoulli beam with point torque boundary conditions described by~\cite{Curtain2012}:
\begin{align*}
\dfrac{\partial^2 y}{\partial t^2} + \dfrac{\partial^4 y}{\partial x^4} - 2 \alpha \dfrac{\partial^3 y}{\partial t \partial x^2} & = u , & & \; \mathrm{in}\;\mathbb{R}_{+}^*\times(0,1) \\
y(t,0) & = 0 , & & \; t \in \mathbb{R}_+^* \\
y(t,1) & = 0 , & & \; t \in \mathbb{R}_+^* \\
\dfrac{\partial^2 y}{\partial x^2} (t,0) & = d_1(t) , & & \; t \in \mathbb{R}_+^* \\
\dfrac{\partial^2 y}{\partial x^2} (t,1) & = d_2(t) , & & \; t \in \mathbb{R}_+^* \\
y(0,x) & = y_0(x) , & & \; x \in (0,1) \\
\dfrac{\partial y}{\partial t} (0,x) & = y_{t0}(x) , & & \; x \in (0,1)
\end{align*}
where $\alpha\in\mathbb{R}_+^*$ is a constant parameter. The functions $u\in\mathcal{C}^0(\mathbb{R}_+;L^2(0,1))$ and $d_1,d_2\in\mathcal{C}^0(\mathbb{R}_+;\mathbb{K})$ are distributed and boundary perturbations, respectively. The functions $y_0 \in H^2(0,1) \cap H_0^1(0,1)$ and $y_{t0} \in L^2(0,1)$ are the initial conditions. 

Introducing the Hilbert space
\begin{equation*}
\mathcal{H} = \left( H^2(0,1) \cap H_0^1(0,1) \right) \times L^2(0,1)
\end{equation*}
with the inner product defined for all $(x_1,x_2), (\hat{x}_1,\hat{x}_2)\in\mathcal{H}$ by
\begin{equation*}
\left<(x_1,x_2) , (\hat{x}_1,\hat{x}_2) \right>_\mathcal{H} = \int_0^1 x_1''(\xi) \overline{\hat{x}_1''(\xi)} + x_2(\xi) \overline{\hat{x}_2(\xi)} \diff \xi ,
\end{equation*}
the distributed parameter system can be written as the abstract system (\ref{def: boundary control system}) with $\mathcal{A}(x_1,x_2) = (x_2,- x_1'''' + 2 \alpha x_2'')$ defined over the domain
\begin{equation*}
D(\mathcal{A}) = \left( H^4(0,1) \cap H_0^1(0,1) \right) \times \left( H^2(0,1) \cap H_0^1(0,1) \right) ,
\end{equation*}
the boundary operator $\mathcal{B}(x_1,x_2) = ( x_1''(0) , x_1''(1) )$ defined over the domain $D(\mathcal{B}) = D(\mathcal{A})$, the state vector $X(t) = (y(t,\cdot),y_t(t,\cdot)) \in \mathcal{H}$, the initial condition $X_0 = (y_0,y_{t0}) \in \mathcal{H}$, $d=(d_1,d_2)\in\mathcal{C}^0(\mathbb{R}_+;\mathbb{K}^2)$, and $U=(0,u)\in\mathcal{C}^0(\mathbb{R}_+;\mathcal{H})$. 

The linear operator $B$ defined such that for all $d=(d_1,d_2)\in\mathbb{K}^2$ and for all $x \in [0,1]$,
\begin{equation*}
(Bd)(x) = \left( \dfrac{d_2-d_1}{6}x^3 + \dfrac{d_1}{2}x^2 - \dfrac{2 d_1 + d_2}{6} x , 0 \right)  ,
\end{equation*}
is a lifting operator associated with $(\mathcal{A},\mathcal{B})$. Following~\cite[Exercise 2.23]{Curtain2012}, it can be shown for $\alpha\in\mathbb{R}_+^*\backslash\{1\}$ that the disturbance-free operator $\mathcal{A}_0$ is a Riesz-spectral operator generating a $C_0$-semigroup of contractions. Thus, $(\mathcal{A},\mathcal{B})$ is a Riesz spectral boundary control system. Furthermore, the eigenvalues of $\mathcal{A}_0$ are given by $\{ - n^2 \pi^2 (\alpha \pm i \sqrt{1-\alpha^2}) , n\in\mathbb{N}^*\}$ when $\alpha\in(0,1)$ while given by $\{ - n^2 \pi^2 (\alpha \pm \sqrt{\alpha^2-1}) , n\in\mathbb{N}^*\}$ when $\alpha > 1$. In both cases, the eigenvalue constraints (\ref{eq: problem setting - confinement conditions}) are satisfied with $\omega_0 = -\alpha\pi^2$ and $\zeta = \alpha^{-1}$ when $\alpha \in (0,1)$ while $\omega_0 = -(\alpha - \sqrt{\alpha^2-1})\pi^2$ and $\zeta = 1$ when $\alpha > 1$.

\subsection{Application of the main results}\label{subsec: ISS estimate derived from the main results}
The adjoint operator $\mathcal{A}_0^*$ is defined over the domain $D(\mathcal{A}_0^*) = D(\mathcal{A}_0)$ by
\begin{equation*}
\mathcal{A}_0^* (x_1,x_2) =
\left(
-x_2 , x_1''''+2 \alpha x_2'' 
\right) .
\end{equation*}
Thus, for initial conditions $y_0 \in H^2(0,1) \cap H_0^1(0,1)$ and $y_{t0} \in L^2(0,1)$, and disturbances $d = (d_1,d_2) \in\mathcal{C}^0(\mathbb{R}_+;\mathbb{K}^2)$ and $u\in\mathcal{C}^0(\mathbb{R}_+;L^2(0,1))$, $X = (x_1,x_2) \in \mathcal{C}^0(\mathbb{R}_+;\mathcal{H})$ is the weak solution of the abstract boundary control system (\ref{def: boundary control system}) associated with $((y_0,y_{t0}),d,(0,u))$ if for all $T>0$ and for all test function $z=(z_1,z_2) \in \mathcal{C}^0([0,T];D(\mathcal{A}^*_0)) \cap \mathcal{C}^1([0,T];\mathcal{H})$ such that $\mathcal{A}_0^* z \in \mathcal{C}^0([0,T];\mathcal{H})$ and $z(T)=0$, the following equality is satisfied: 
\begin{align*}
& \int_0^T \int_0^1 x_1''(t,\xi) \overline{ \left\{ - z_2'' + \left[ \dfrac{\mathrm{d} z_1}{\mathrm{d}t} \right]'' \right\} (t,\xi) } \diff\xi \diff t \\
& + \int_0^T \int_0^1 x_2(t,\xi) \overline{ \left\{ z_1'''' + 2 \alpha z_2'' + \dfrac{\mathrm{d} z_2}{\mathrm{d}t} \right\} (t,\xi) } \diff\xi \diff t \\ 
& = - \int_0^1 y_0''(\xi) \overline{z_1''(0,\xi)} \diff\xi
- \int_0^1 y_{t0}(\xi) \overline{z_2(0,\xi)} \diff\xi \\
& \phantom{=}\; + \int_0^T d_1(t) \overline{z_2'(t,0)} \diff t  - \int_0^T d_2(t) \overline{z_2'(t,1)} \diff t \\
& \phantom{=}\; - \int_0^T \int_0^1 u(t,\xi) \overline{z_2(t,\xi)} \diff\xi \diff t .
\end{align*}

As the chosen lifting operator satisfies $\mathcal{A}B = 0$, the constants $C_i$ provided by both Theorems~\ref{th: ISS classical solutions} and~\ref{thm: constant C1 independant of B} are identical. Following~\cite[Exercise 2.23]{Curtain2012} and using the same approach that the one used in~\cite{lhachemi2018input} for establishing the constants $m_R$ and $M_R$ related to the Riesz basis, we have:
\begin{enumerate}
\item Case $\alpha \in (0,1)$: $\kappa_0 = \alpha\pi^2$, $\zeta = \alpha^{-1}$, $m_R = 1 - \alpha$, $M_R = 1 + \alpha$;
\item Case $\alpha > 1$: $\kappa_0 = (\alpha - \sqrt{\alpha^2-1})\pi^2$, $\zeta = 1$, $m_R = 1 - \alpha^{-1}$, $M_R = 1 + \alpha^{-1}$.
\end{enumerate}
The boundary disturbance evolves into the two dimensional ($m=2$) space $(\mathbb{K}^2,\Vert\cdot\Vert_2)$ endowed with the usual euclidean norm. By selecting $\mathcal{E}=\{e_1,e_2\}$ as the canonical basis of $\mathbb{K}^2$, we obtain $c(\mathcal{E})=1$ and $\Vert Be_1 \Vert_\mathcal{H} = \Vert Be_2 \Vert_\mathcal{H} = 1/\sqrt{3}$. Thus, the constants of the ISS estimate are given by
\begin{equation}\label{eq: ISS constants V1 alpha>1}
C_0 = \sqrt{\dfrac{1+\alpha}{1-\alpha}} , \;\; C_1 = \dfrac{2}{\alpha\sqrt{3}}\sqrt{\dfrac{1+\alpha}{1-\alpha}} , \;\; C_2 = \dfrac{1}{\alpha\pi^2}\sqrt{\dfrac{1+\alpha}{1-\alpha}}   
\end{equation}
when $\alpha \in (0,1)$, while
\begin{equation}\label{eq: ISS constants V1 alpha<1}
C_0 = \sqrt{\dfrac{\alpha+1}{\alpha-1}} , \;\; C_1 = \dfrac{2}{\sqrt{3}}\sqrt{\dfrac{\alpha+1}{\alpha-1}} , 
\end{equation}
\begin{equation*}
C_2 = \dfrac{1}{(\alpha - \sqrt{\alpha^2-1})\pi^2}\sqrt{\dfrac{\alpha+1}{\alpha-1}}  
\end{equation*}
when $\alpha > 1$. 

We observe that the ISS constants provided by Theorems~\ref{th: ISS classical solutions} and~\ref{thm: constant C1 independant of B} diverge to $+ \infty$ when the quantity $\alpha \rightarrow 1$. This is due to the fact that the disturbance-free operator $\mathcal{A}_0$ losses its Riesz-spectral properties when $\alpha = 1$ with in particular $m_R \underset{\alpha \rightarrow 1}{\longrightarrow} 0$. Thus, the obtained ISS constants become conservative when $\alpha \rightarrow 1$ due to the fact that the constants of the Riesz basis are such that $M_R / m_R \underset{\alpha \rightarrow 1}{\longrightarrow} +\infty$. An approach for reducing such a conservatism is discussed in the next subsection.

It is interesting to evaluate the impact of an increased damping term $\alpha$. Constants $C_0$ and $C_1$ are decreasing functions of $\alpha > 1$ and we have $C_0 \underset{\alpha \rightarrow +\infty}{\longrightarrow} 1$ and $C_1 \underset{\alpha \rightarrow +\infty}{\longrightarrow} 2/\sqrt{3}$. We also observe that $C_2 \underset{\alpha \rightarrow +\infty}{\longrightarrow} +\infty$; this divergent behavior is induced by the fact that Theorems~\ref{th: ISS classical solutions} and~\ref{thm: constant C1 independant of B} deal with distributed perturbations evolving in the full space $\mathcal{H}$. However, the distributed disturbance evolves in the subspace $\{0\} \times L^2(0,1)$. To provide a tighter constant $C_2$, we directly estimate (\ref{eq: trajectory - integral term U}) with the sparse structure $U=(0,u)$. This yields the following version of the ISS constant:
\begin{equation*}
C_2 = 
\sqrt{ M_R \sum\limits_{\substack{n \in \mathbb{N}^* \\ \epsilon \in \{-1,+1\}}} \dfrac{\Vert \psi_{n,\epsilon}^{2} \Vert_{L^2(0,1)}^2}{\left| \operatorname{Re} \lambda_{n,\epsilon} \right|^2} }
\end{equation*}
where $\psi_{n,\varepsilon}^2$ denotes the second component of $\psi_{n,\varepsilon}$.
For $\alpha>1$, the eigenvalues of the disturbance-free operator $\mathcal{A}_0$ are given by $\lambda_{n,\epsilon} = -n^2\pi^2(\alpha+\epsilon\sqrt{\alpha^2-1}) \in\mathbb{R}_-^*$, $n\in\mathbb{N}^*$ and $\varepsilon\in\{-1,+1\}$. The corresponding eigenvectors are given by
\begin{equation*}
\phi_{n,\varepsilon} = \dfrac{1}{n^2 \pi^2 \sqrt{\alpha \left( \alpha + \varepsilon \sqrt{\alpha^2-1} \right)}} 
\begin{pmatrix}
\sin(n \pi \cdot) \\ \lambda_{n,\varepsilon} \sin(n \pi \cdot)
\end{pmatrix} .
\end{equation*}
The associated biorthogonal sequence is given by
\begin{equation*}
\psi_{n,\varepsilon} = \dfrac{2 \sqrt{\alpha \left(\alpha+\varepsilon\sqrt{\alpha^2-1}\right)}}{n^2 \pi^2 \left( 1 - (\alpha+\epsilon\sqrt{\alpha^2-1})^2 \right)} 
\begin{pmatrix}
\sin(n \pi \cdot) \\ - \lambda_{n,\varepsilon} \sin(n \pi \cdot)
\end{pmatrix} .
\end{equation*}
Now, straightforward computations show that\footnote{We used $\sum\limits_{n\geq 1} \dfrac{1}{n^4} = \dfrac{\pi^4}{90}$.} 
\begin{equation*}
C_2 = \dfrac{1}{3 \sqrt{10}} \sqrt{\dfrac{\alpha}{\alpha-1}} 
\underset{\alpha \rightarrow +\infty}{\longrightarrow} \dfrac{1}{3 \sqrt{10}} .
\end{equation*}

\subsection{Improvement in the neighborhood of $\alpha = 1$ by adaptation of the spectral decomposition}\label{subsec: Improvement in the neighborhood of alpha = 1}
\begin{figure}
\centering
\includegraphics[width=3.3in,height=1in]{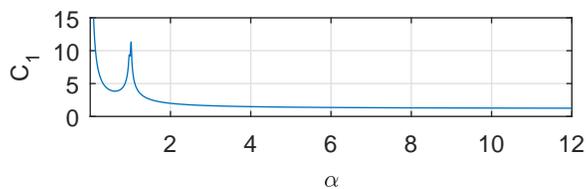}
\caption{ISS constant $C_1$ for the Euler-Bernoulli beam}
\label{fig: beam ISS constant C1}
\end{figure}

We propose to adapt the spectral decomposition (\ref{eq: diff equation c_n}) in order to reduce the conservatism observed in the neighborhood of $\alpha = 1$ for the constants $C_0$ and $C_1$. To do so, we introduce $\psi_n = \dfrac{1}{n^2 \pi^2} \left( \sin(n\pi\cdot) , n^2 \pi^2 \sin(n\pi\cdot) \right) \in D(\mathcal{A}_0^*)$, $n \in \mathbb{N}^*$, which are the eigenvectors of $\mathcal{A}_0^*$ in the case $\alpha = 1$ associated with the eigenvalues $\mu_n = -n^2 \pi^2$. As $\left\{ \psi_n , \; n \in \mathbb{N} \right\}$ is not maximal in $\mathcal{H}$, we introduce the vectors $\psi_n^d = \dfrac{1}{n^2 \pi^2} \left( \sin(n\pi\cdot) , - n^2 \pi^2 \sin(n\pi\cdot) \right) \in D(\mathcal{A}_0^*)$. Straightforward computations show that $\left\{ \psi_n ,\psi_n^d , \; n \in \mathbb{N}^* \right\}$ is a Hilbert basis of $\mathcal{H}$ and, for all $n \in \mathbb{N}^*$,
\begin{align*}
\mathcal{A}_0^* \psi_n & = - \alpha n^2 \pi^2 \psi_n + (\alpha-1) n^2 \pi^2 \psi_n^d , \\
\mathcal{A}_0^* \psi_n^d & = (\alpha+1) n^2 \pi^2 \psi_n - \alpha n^2 \pi^2 \psi_n^d .
\end{align*}
Then, introducing $c_n(t) \triangleq \left< X(t) , \psi_n \right>_\mathcal{H}$ and $c_n^d(t) \triangleq \left< X(t) , \psi_n^d \right>_\mathcal{H}$, we have $\Vert X(t) \Vert_\mathcal{H}^2 = \sum\limits_{n \geq 1} \vert c_n(t) \vert^2 + \vert c_n^d(t) \vert^2$. Furthermore, by using the same approach as the one employed to derive (\ref{eq: diff equation c_n}), the following decomposition holds true for any classical solution $X$ associated with an initial condition $X_0$ and a boundary disturbance $d$, 
\begin{equation*}
\begin{bmatrix} \dot{c}_n(t) \\ \dot{c}_n^d(t) \end{bmatrix} 
= n^2 \pi^2 M \begin{bmatrix} c_n(t) \\ c_n^d(t) \end{bmatrix} 
- n^2 \pi^2 M \begin{bmatrix} \left< D(t) , \psi_n \right>_\mathcal{H} \\ \left< D(t) , \psi_n^d \right>_\mathcal{H} \end{bmatrix}
\end{equation*} 
with
\begin{equation*}
M \triangleq 
\begin{bmatrix}
- \alpha & \alpha-1 \\ \alpha+1 & -\alpha
\end{bmatrix}
\end{equation*}
and where $D(t)$ is defined by (\ref{eq: project boundary disturbance over stationnary sol}). Then, for all $t \geq 0$, we have
\begin{align*}
\begin{bmatrix} c_n(t) \\ c_n^d(t) \end{bmatrix} 
& = e^{n^2 \pi^2 M t} \begin{bmatrix} c_n(0) \\ c_n^d(0) \end{bmatrix} \\
& \phantom{=}\, - n^2 \pi^2 M \int_0^t e^{n^2 \pi^2 M (t-\tau)} \begin{bmatrix} \left< D(\tau) , \psi_n \right>_\mathcal{H} \\ \left< D(\tau) , \psi_n^d \right>_\mathcal{H} \end{bmatrix} \diff\tau .
\end{align*} 
In the case $\alpha > 1$, we have
\begin{equation*}
e^{n^2 \pi^2 M t} 
= e^{- \alpha n^2 \pi^2 t}
\begin{bmatrix}
\cosh(\theta_{n}t) & \sqrt{\dfrac{\alpha-1}{\alpha+1}} \sinh(\theta_{n}t) \\
\sqrt{\dfrac{\alpha+1}{\alpha-1}} \sinh(\theta_{n}t) & \cosh(\theta_{n}t)
\end{bmatrix}
\end{equation*}
where $\theta_{n} \triangleq n^2\pi^2\sqrt{\alpha^2-1}$. Then, with a similar approach as the one employed in the proof of Theorem~\ref{th: ISS classical solutions} and using the inequalities $\cosh(x) \leq e^x$, $\sinh(x) \leq e^{x}/2$, $\sinh(x) \leq x \cosh(x)$, and $x e^{-\gamma x} \leq e^{-1}/\gamma$ that hold for all $x \geq 0$ and $\gamma > 0$, we obtain for any $\epsilon \in (0,1)$ the following versions of the ISS constants:
\begin{align}
\kappa_0 & = (1-\epsilon) (\alpha-\sqrt{\alpha^2-1})\pi^2 , \nonumber \\
C_0 & = 2 \max\left( 1 , \dfrac{e^{-1}(\alpha+1)}{\epsilon(\alpha - \sqrt{\alpha^2-1})} \right) , \label{eq: ISS constants V2 alpha>1} \\
C_1 & = \dfrac{4(\alpha + 1)(2\alpha-\sqrt{\alpha^2-1})}{\sqrt{3}(\alpha-\sqrt{\alpha^2-1})^2} . \nonumber
\end{align}
Similar computations show that the above constants are also valid in the case $\alpha = 1$ while the case $\alpha \in (0,1)$ yields the following versions of the constants:
\begin{equation*}
\kappa_0 = (1-\epsilon) \alpha\pi^2 ,
\end{equation*}
\begin{equation}\label{eq: ISS constants V2 alpha<1}
C_0 = 2 \max\left( 1 , \dfrac{e^{-1}(\alpha+1)}{\epsilon\alpha} \right) ,
\quad
C_1 = \dfrac{8(\alpha + 1)}{\sqrt{3}\alpha}.
\end{equation}
Comparing to (\ref{eq: ISS constants V1 alpha>1}-\ref{eq: ISS constants V1 alpha<1}), the above ISS constants (\ref{eq: ISS constants V2 alpha>1}-\ref{eq: ISS constants V2 alpha<1}) are less conservative in a small neighborhood of $\alpha = 1$ because finite in $\alpha = 1$. However, outside such a neighborhood, they become more conservative. For instance, the above version of $C_1$ given by (\ref{eq: ISS constants V2 alpha>1}-\ref{eq: ISS constants V2 alpha<1}), which is such that $C_1 \underset{\alpha \rightarrow +\infty}{\sim} 16 \alpha^4 / \sqrt{3}$, is only less conservative than the version (\ref{eq: ISS constants V1 alpha>1}-\ref{eq: ISS constants V1 alpha<1}) over the (approximate) range $\alpha \in (0.967,1.020)$. The evolution of the constant $C_1$ obtained by taking the minimum of the versions (\ref{eq: ISS constants V1 alpha>1}-\ref{eq: ISS constants V1 alpha<1}) and (\ref{eq: ISS constants V2 alpha>1}-\ref{eq: ISS constants V2 alpha<1}) is depicted in Fig.~\ref{fig: beam ISS constant C1}.  

\section{Conlusion}\label{sec: conclusion}
This paper discussed the establishment of Input-to-State Stability (ISS) estimates for a class of Riesz-spectral boundary control system with respect to both boundary and distributed perturbations. First, a spectral decomposition depending only on the boundary disturbance (but not on its time derivative) was obtained by projecting the system trajectories over an adequate Riesz basis. This was used to derive an ISS estimate with respect to classical solutions. Then, in order to relax the regularity assumptions required for assessing the existence of classical solutions, a concept of weak solution that applies to a large class of boundary control systems (which is not limited to Riesz-spectral ones) has been introduced under a variational formulation. Assuming that an ISS estimate holds true with respect to classical solutions, various properties of the weak solutions were derived, including their existence and uniqueness, as well as their ISS property.


\bibliographystyle{plain}        
\bibliography{autosam}           

\begin{thebibliography}{10}

\bibitem{argomedo2013strict}
Federico~Bribiesca Argomedo, Christophe Prieur, Emmanuel Witrant, and Sylvain
  Br{\'e}mond.
\newblock A strict control {L}yapunov function for a diffusion equation with
  time-varying distributed coefficients.
\newblock {\em IEEE Transactions on Automatic Control}, 58(2):290--303, 2013.

\bibitem{argomedo2012d}
Federico~Bribiesca Argomedo, Emmanuel Witrant, and Christophe Prieur.
\newblock {$D^1$}-{I}nput-to-state stability of a time-varying nonhomogeneous
  diffusive equation subject to boundary disturbances.
\newblock In {\em American Control Conference (ACC), 2012}, pages 2978--2983.
  IEEE, 2012.

\bibitem{ball1977strongly}
JM~Ball.
\newblock Strongly continuous semigroups, weak solutions, and the variation of
  constants formula.
\newblock {\em Proceedings of the American Mathematical Society},
  63(2):370--373, 1977.

\bibitem{brezis2010functional}
Haim Brezis.
\newblock {\em Functional Analysis, Sobolev Spaces and Partial Differential
  Equations}.
\newblock Springer Science \& Business Media, 2010.

\bibitem{christensen2016introduction}
Ole Christensen et~al.
\newblock {\em An Introduction to Frames and {Riesz} Bases}.
\newblock Springer, 2016.

\bibitem{coron2004global}
Jean-Michel Coron and Emmanuel Tr{\'e}lat.
\newblock Global steady-state controllability of one-dimensional semilinear
  heat equations.
\newblock {\em SIAM Journal on Control and Optimization}, 43(2):549--569, 2004.

\bibitem{coron2006global}
Jean-Michel Coron and Emmanuel Tr{\'e}lat.
\newblock Global steady-state stabilization and controllability of {1D}
  semilinear wave equations.
\newblock {\em Communications in Contemporary Mathematics}, 8(04):535--567,
  2006.

\bibitem{Curtain2012}
R.~F. Curtain and H.~Zwart.
\newblock {\em An Introduction to Infinite-Dimensional Linear Systems Theory},
  volume~21.
\newblock Springer Science \& Business Media, 2012.

\bibitem{dashkovskiy2011local}
Sergey Dashkovskiy and Andrii Mironchenko.
\newblock Local {ISS} of reaction-diffusion systems.
\newblock {\em IFAC Proceedings Volumes}, 44(1):11018--11023, 2011.

\bibitem{dashkovskiy2013input}
Sergey Dashkovskiy and Andrii Mironchenko.
\newblock Input-to-state stability of infinite-dimensional control systems.
\newblock {\em Mathematics of Control, Signals, and Systems}, 25(1):1--35,
  2013.

\bibitem{emirsjlow2000pdes}
Zbigniew Emirsjlow and Stuart Townley.
\newblock From {PDEs} with boundary control to the abstract state equation with
  an unbounded input operator: a tutorial.
\newblock {\em European Journal of Control}, 6(1):27--49, 2000.

\bibitem{endo2017boundary}
Takahiro Endo, Fumitoshi Matsuno, and Yingmin Jia.
\newblock Boundary cooperative control by flexible {Timoshenko} arms.
\newblock {\em Automatica}, 81:377--389, 2017.

\bibitem{henikl2016infinite}
Johannes Henikl, Wolfgang Kemmetm{\"u}ller, Thomas Meurer, and Andreas Kugi.
\newblock Infinite-dimensional decentralized damping control of large-scale
  manipulators with hydraulic actuation.
\newblock {\em automatica}, 63:101--115, 2016.

\bibitem{jacob2016input}
Birgit Jacob, Robert Nabiullin, Jonathan Partington, and Felix Schwenninger.
\newblock On input-to-state-stability and integral input-to-state-stability for
  parabolic boundary control systems.
\newblock In {\em 2016 IEEE 55th Conference on Decision and Control (CDC)},
  pages 2265--2269. IEEE, 2016.

\bibitem{jacob2018infinite}
Birgit Jacob, Robert Nabiullin, Jonathan~R Partington, and Felix~L
  Schwenninger.
\newblock Infinite-dimensional input-to-state stability and orlicz spaces.
\newblock {\em SIAM Journal on Control and Optimization}, 56(2):868--889, 2018.

\bibitem{jacob2018continuity}
Birgit Jacob, Felix~L Schwenninger, and Hans Zwart.
\newblock On continuity of solutions for parabolic control systems and
  input-to-state stability.
\newblock {\em Journal of differential equations, in press}, 2018.

\bibitem{karafyllis2018boundary}
Iasson Karafyllis, Maria Kontorinaki, and Miroslav Krstic.
\newblock Boundary-to-displacement asymptotic gains for wave systems with
  {Kelvin-Voigt} damping.
\newblock {\em arXiv preprint arXiv:1807.06549}, 2018.

\bibitem{karafyllis2016input}
Iasson Karafyllis and Miroslav Krstic.
\newblock Input-to state stability with respect to boundary disturbances for
  the {1-D} heat equation.
\newblock In {\em Decision and Control (CDC), 2016 IEEE 55th Conference on},
  pages 2247--2252. IEEE, 2016.

\bibitem{karafyllis2016iss}
Iasson Karafyllis and Miroslav Krstic.
\newblock {ISS} with respect to boundary disturbances for {1-D} parabolic
  {PDEs}.
\newblock {\em IEEE Transactions on Automatic Control}, 61(12):3712--3724,
  2016.

\bibitem{karafyllis2017iss}
Iasson Karafyllis and Miroslav Krstic.
\newblock {ISS} in different norms for {1-D} parabolic {PDEs} with boundary
  disturbances.
\newblock {\em SIAM Journal on Control and Optimization}, 55(3):1716--1751,
  2017.

\bibitem{karafyllis2019preview}
Iasson Karafyllis and Miroslav Krstic.
\newblock {\em Input-to-State Stability for {PDEs}}.
\newblock Springer, 2019.

\bibitem{lhachemi2018boundaryAutomatica}
Hugo Lhachemi, David Saussi{\'e}, and Guchuan Zhu.
\newblock Boundary feedback stabilization of a flexible wing model under
  unsteady aerodynamic loads.
\newblock {\em Automatica}, 97:73--81, 2018.

\bibitem{lhachemi2018input}
Hugo Lhachemi, David Saussi{\'e}, Guchuan Zhu, and Robert Shorten.
\newblock Input-to-state stability of a clamped-free damped string in the
  presence of distributed and boundary disturbances.
\newblock {\em arXiv preprint arXiv:1807.11696}, 2018.

\bibitem{mazenc2011strict}
Fr{\'e}d{\'e}ric Mazenc and Christophe Prieur.
\newblock Strict {L}yapunov functions for semilinear parabolic partial
  differential equations.
\newblock {\em Mathematical Control and Related Fields}, 1(2):231--250, 2011.

\bibitem{mironchenko2016local}
Andrii Mironchenko.
\newblock Local input-to-state stability: Characterizations and
  counterexamples.
\newblock {\em Systems \& Control Letters}, 87:23--28, 2016.

\bibitem{mironchenko2014integral}
Andrii Mironchenko and Hiroshi Ito.
\newblock Integral input-to-state stability of bilinear infinite-dimensional
  systems.
\newblock In {\em Decision and Control (CDC), 2014 IEEE 53rd Annual Conference
  on}, pages 3155--3160. IEEE, 2014.

\bibitem{mironchenko2015construction}
Andrii Mironchenko and Hiroshi Ito.
\newblock Construction of {L}yapunov functions for interconnected parabolic
  systems: an {iISS} approach.
\newblock {\em SIAM Journal on Control and Optimization}, 53(6):3364--3382,
  2015.

\bibitem{mironchenko2019monotonicity}
Andrii Mironchenko, Iasson Karafyllis, and Miroslav Krstic.
\newblock Monotonicity methods for input-to-state stability of nonlinear
  parabolic pdes with boundary disturbances.
\newblock {\em SIAM Journal on Control and Optimization}, 57(1):510--532, 2019.

\bibitem{mironchenko2016restatements}
Andrii Mironchenko and Fabian Wirth.
\newblock Restatements of input-to-state stability in infinite dimensions: what
  goes wrong.
\newblock In {\em Proc. of 22th International Symposium on Mathematical Theory
  of Systems and Networks (MTNS 2016)}, pages 667--674, 2016.

\bibitem{mironchenko2017characterizations}
Andrii Mironchenko and Fabian Wirth.
\newblock Characterizations of input-to-state stability for
  infinite-dimensional systems.
\newblock {\em IEEE Transactions on Automatic Control}, 63(6):1692--1707, 2018.

\bibitem{prieur2012iss}
Christophe Prieur and Fr{\'e}d{\'e}ric Mazenc.
\newblock {ISS}-{L}yapunov functions for time-varying hyperbolic systems of
  balance laws.
\newblock {\em Mathematics of Control, Signals, and Systems}, 24(1-2):111--134,
  2012.

\bibitem{schmid2018stabilization}
Jochen Schmid and Hans Zwart.
\newblock Stabilization of port-{H}amiltonian systems by nonlinear boundary
  control in the presence of disturbances.
\newblock {\em arXiv preprint arXiv:1804.10598}, 2018.

\bibitem{sontag1989smooth}
Eduardo~D Sontag.
\newblock Smooth stabilization implies coprime factorization.
\newblock {\em IEEE transactions on automatic control}, 34(4):435--443, 1989.

\bibitem{tanwani2017disturbance}
Aneel Tanwani, Christophe Prieur, and Sophie Tarbouriech.
\newblock Disturbance-to-state stabilization and quantized control for linear
  hyperbolic systems.
\newblock {\em arXiv preprint arXiv:1703.00302}, 2017.

\bibitem{tucsnak2009observation}
Marius Tucsnak and George Weiss.
\newblock {\em Observation and Control for Operator Semigroups}.
\newblock Springer Science \& Business Media, 2009.

\bibitem{zheng2017giorgi}
Jun Zheng and Guchuan Zhu.
\newblock A {De Giorgi} iteration-based approach for the establishment of {ISS}
  properties for {B}urgers' equation with boundary and in-domain disturbances.
\newblock {\em IEEE Transactions on Automatic Control, in press}, 2018.

\bibitem{zheng2017input}
Jun Zheng and Guchuan Zhu.
\newblock Input-to-state stability with respect to boundary disturbances for a
  class of semi-linear parabolic equations.
\newblock {\em Automatica}, 97:271--277, 2018.

\end{thebibliography}



\appendix

\section{Proof of Lemma~\ref{lemma: weak solution uniqueness}}\label{annex: proof lemma uniqueness weak solutions}    

By linearity, we must show that if $X \in \mathcal{C}^0(\mathbb{R}_+;\mathcal{H})$ satisfies
\begin{equation}\label{eq: weak solution kernel}
\int_0^T \left< X(t), \mathcal{A}_0^* z(t) + \dfrac{\mathrm{d} z}{\mathrm{d} t}(t) \right>_\mathcal{H} \diff t = 0 
\end{equation}
for all $T>0$ and for all $z \in \mathcal{C}^0([0,T];D(\mathcal{A}^*_0)) \cap \mathcal{C}^1([0,T];\mathcal{H})$ such that $\mathcal{A}_0^* z \in \mathcal{C}^0([0,T];\mathcal{H})$ and $z(T)=0$, then $X=0$. Denoting by $S$ the $C_0$-semigroup generated by $\mathcal{A}_0$, then $S^*$ is the $C_0$-semigroup generated by $\mathcal{A}_0^*$ (see, e.g., ~\cite[Thm 2.2.6]{Curtain2012}). Let $z_0 \in D(\mathcal{A}_0^*)$ and $\alpha > 0$ be arbitrarily given. For any given $T>0$, we consider the function $z_{z_0,\alpha,T}$ defined for any $t \in [0,T]$ by $z_{z_0,\alpha,T}(t) = S^*(\alpha t) z_0 - S^*(\alpha T) z_0$. As $z_0 \in D(\mathcal{A}_0^*)$, we obtain that for any $t \geq 0$, $S^*(\alpha t) z_0 \in D(\mathcal{A}_0^*)$, $\mathcal{A}_0^* S^*(\alpha t) z_0 = S^*(\alpha t) \mathcal{A}_0^* z_0$, and 
\begin{equation*}
\dfrac{\mathrm{d} z_{z_0,\alpha,T}}{\mathrm{d} t}(t) = \alpha \mathcal{A}_0^* S^*(\alpha t) z_0 = \alpha S^*(\alpha t) \mathcal{A}_0^* z_0 .
\end{equation*}
Thus, $z_{z_0,\alpha,T}$ is a test function over $[0,T]$ and (\ref{eq: weak solution kernel}) shows:
\begin{align*}
& (\alpha+1) \int_0^T \left< X(t), S^*(\alpha t) \mathcal{A}_0^* z_0\right>_\mathcal{H} \diff t \\
& \phantom{=======} = \int_0^T \left< X(t), S^*(\alpha T) \mathcal{A}_0^* z_0\right>_\mathcal{H} \diff t .
\end{align*}
Using the definition of the adjoint operator, the fact that $S(\alpha T) \in \mathcal{L}(\mathcal{H})$, and the properties of the Bochner integral, the equation above is equivalent to
\begin{align*}
& \left< (\alpha+1) \int_0^T  S(\alpha t) X(t) \diff t - S(\alpha T) \int_0^T X(t) \diff t , \mathcal{A}_0^* z_0\right>_\mathcal{H} \\
& = 0.
\end{align*}
Because $\mathcal{H}$ is a Hilbert space with $\mathcal{A}_0$ closed and densely defined, we have that $\mathrm{ker}\left( \mathcal{A}_0 \right)^\bot = \overline{\mathrm{R}\left( \mathcal{A}_0^* \right)}$ (see, e.g., \cite[Chap.~2, Rem.~17]{brezis2010functional}). Since $\mathcal{A}_0$ is assumed to be injective, this yields that $\overline{\mathrm{R}\left( \mathcal{A}_0^* \right)} = \mathcal{H}$. Consequently, we obtain that the following equality holds true for all $\alpha > 0$ and $T \geq 0$
\begin{equation}\label{eq: uniqueness weak solutions - interm eq}
(\alpha+1) \int_0^T  S(\alpha t) X(t) \diff t = S(\alpha T) \int_0^T X(t) \diff t .
\end{equation}
This implies that, for any $h > 0$,
\begin{align*}
& \dfrac{\alpha + 1}{h} \int_{T}^{T+h} S(\alpha t) X(t) \diff t \\
& \phantom{===} = S(\alpha(T+h)) \left\{ \dfrac{1}{h} \int_{T}^{T+h} X(t) \diff t \right\} \\
& \phantom{===} \phantom{=}\, + \dfrac{S(\alpha h) - I_\mathcal{H}}{\alpha h} \left\{ \alpha S(\alpha T) \int_0^T X(t) \diff t \right\} . 
\end{align*}
As $t \rightarrow X(t)$ and $t \rightarrow S(\alpha t)X(t)$ are continuous over $\mathbb{R}_+$, we obtain by the continuity property of the $C_0$-semigroups that
\begin{align*}
& \underset{h \rightarrow 0^+}{\lim}  \dfrac{S(\alpha h) - I_\mathcal{H}}{\alpha h} \left\{ S(\alpha T) \int_0^T X(t) \diff t \right\} \\
& = \dfrac{\alpha + 1}{\alpha} S(\alpha T) X(T) - \dfrac{1}{\alpha} S(\alpha T) X(T) \\
& = S(\alpha T) X(T) .
\end{align*}
Thus we have $S(\alpha T) \int_0^T X(t) \diff t \in D(\mathcal{A}_0)$ and
\begin{equation*}
\mathcal{A}_0 \left\{ S(\alpha T) \int_0^T X(t) \diff t \right\}
= S(\alpha T) X(T).
\end{equation*}
From (\ref{eq: uniqueness weak solutions - interm eq}), we deduce that $\int_0^T  S(\alpha t) X(t) \diff t \in D(\mathcal{A}_0)$ and
\begin{equation*}
(\alpha + 1) \mathcal{A}_0 \int_0^T  S(\alpha t) X(t) \diff t = S(\alpha T)X(T)
\end{equation*}
for all $\alpha > 0$ and $T \geq 0$. Introducing $y_\alpha(T) = \int_0^T  S(\alpha t) X(t) \diff t$ and noting that $t \rightarrow S(\alpha t) X(t)$ is continuous over $\mathbb{R}_+$, we obtain that $y_\alpha$ satisfies over $\mathbb{R}_+$ the differential equation $\dfrac{\mathrm{d}y_\alpha}{\mathrm{d}T} = (\alpha+1)\mathcal{A}_0 y_\alpha$ with the initial condition $y_\alpha(0) = 0$. As $\mathcal{A}_0$ generates the $C_0$-semigroup $S$ and $\alpha + 1 > 0$, we deduce that $(\alpha+1)\mathcal{A}_0$ generates the $C_0$-semigroup $S((\alpha+1)\cdot)$. Thus we have $y_\alpha = S((\alpha+1)\cdot) y_\alpha(0) = 0$. By taking the time derivative of $y_\alpha$, we deduce that $S(\alpha t)X(t) = 0$ for all $\alpha > 0$ and $t \geq 0$. From the continuity property of the $C_0$-semigroups, we obtain by letting $\alpha \rightarrow 0^+$ that $X(t)=0$ for all $t \geq 0$.\qed

\section{Proof of Theorem~\ref{th: extended main result}}\label{annex: proof thm extended main result}    

To establish the uniqueness part, we only need to show that $\mathcal{A}_0$ is injective. In that case, the conclusion will follow from the application of Lemma~\ref{lemma: weak solution uniqueness}. Let $x_0 \in \mathrm{ker}(\mathcal{A}_0)$ be arbitrarily given. Introducing $X(t) = x_0$ for all $t \geq 0$, $X_0 = x_0$, $d=0$, and $U = 0$, one has for all $t \geq 0$, $\dfrac{\mathrm{d}X}{\mathrm{d}t}(t) = 0 = \mathcal{A}_0 x_0 = \mathcal{A} X(t) + U(t)$, $\mathcal{B}X(t) = \mathcal{B}x_0 = 0 = d(t)$, and $X(0)=x_0$. Thus $X$ is the classical solution associated with $(X_0,d,U)=(x_0,0,0)$. The ISS estimate (\ref{eq: th weak solutions - ISS}) gives  $\left\Vert x_0 \right\Vert_\mathcal{H} \leq \beta(\left\Vert x_0 \right\Vert_\mathcal{H},t) 
\underset{t \rightarrow +\infty}{\longrightarrow} 0$. This yields $x_0=0$, ensuring the injectivity of $\mathcal{A}_0$.

To show the existence part, let an initial condition $X_0 \in \mathcal{H}$, and disturbances $d\in\mathcal{C}^0(\mathbb{R}_+;\mathbb{K}^m)$ and $U\in\mathcal{C}^0(\mathbb{R}_+;\mathcal{H})$ be arbitrarily given. We also consider an arbitrarily given lifting operator $B$ associated with $(\mathcal{A},\mathcal{B})$.

\textit{Step~1: Construction of a weak solution candidate $X \in \mathcal{C}^0([0,T];\mathcal{H})$ by density arguments.}

Let $T > 0$ be arbitrarily given. As $\mathcal{C}^2([0,T];\mathbb{K}^m)$ and $\mathcal{C}^1([0,T];\mathcal{H})$ are dense in $\mathcal{C}^0([0,T];\mathbb{K}^m)$ and $\mathcal{C}^0([0,T];\mathcal{H})$, respectively, there exist sequences $(d_n)_n \in \mathcal{C}^2([0,T];\mathbb{K}^m)^\mathbb{N}$ and $(U_n)_n \in \mathcal{C}^1([0,T];\mathcal{H})^\mathbb{N}$ (we can actually use here approximating sequences of smooth functions, i.e., of class $\mathcal{C}^\infty$) such that 
\begin{align*}
\Vert d_n - d \Vert_{\mathcal{C}^0([0,T];\mathbb{K}^m)} & \underset{n \rightarrow +\infty}{\longrightarrow} 0 , \\
\Vert U_n - U \Vert_{\mathcal{C}^0([0,T];\mathcal{H})} & \underset{n \rightarrow +\infty}{\longrightarrow} 0 .
\end{align*}
Now, as $\overline{D(\mathcal{A}_0)} = \mathcal{H}$, there exists $\left( \tilde{X}_{0,n} \right)_n \in (D(\mathcal{A}_0))^\mathbb{N}$ such that $\tilde{X}_{0,n} \underset{n \rightarrow +\infty}{\longrightarrow} X_0 - Bd(0)$. Introducing $X_{0,n} \triangleq \tilde{X}_{0,n} + B d_n(0) \in D(\mathcal{A})$, the bounded nature of $B$ gives $X_{0,n} \underset{n \rightarrow +\infty}{\longrightarrow} X_0$. Recalling that $D(\mathcal{A}_0) \subset \mathrm{ker}(\mathcal{B})$ and $\mathcal{B}B = I_\mathbb{K}$, we get $\mathcal{B} X_{0,n} = \mathcal{B} \tilde{X}_{0,n} + \mathcal{B} B d_n(0) = d_n(0)$.

For any $n \in \mathbb{N}$, let $X_n \in C^0([0,T];D(\mathcal{A})) \cap C^1([0,T];\mathcal{H})$ be the classical solution of the abstract system (\ref{def: boundary control system}) over $[0,T]$ associated with $(X_{0,n},d_n,U_n)$. For any $n,m\in\mathbb{N}$, by linearity, $X_n - X_m$ is the unique classical solution of the abstract system (\ref{def: boundary control system}) over $[0,T]$ associated with $(X_{0,n} - X_{0,m},d_n - d_m,U_n - U_m)$. Thus the ISS estimate for classical solutions (\ref{eq: th weak solutions - ISS}) yields for all $n,m\in\mathbb{N}$, 
\begin{align*}
\Vert X_n - X_m \Vert_{\mathcal{C}^0([0,T];\mathcal{H})} \leq
& \beta \left( \Vert X_{0,n} - X_{0,m} \Vert_\mathcal{H},0 \right) \\
& + \gamma_1 \left( \Vert d_n - d_m \Vert_{\mathcal{C}^0([0,T];\mathbb{K}^m)} \right) \\ 
& + \gamma_2 \left( \Vert U_n - U_m \Vert_{\mathcal{C}^0([0,T];\mathcal{H})} \right) .
\end{align*}  
Since $(X_{0,n})_n$, $(d_n)_n$, and $(U_n)_n$ are Cauchy sequences, so is $(X_n)_n$. As $\mathcal{C}^0([0,T];\mathcal{H})$ is a Banach space, there exists $X \in \mathcal{C}^0([0,T];\mathcal{H})$ such that $\Vert X_n - X \Vert_{\mathcal{C}^0([0,T];\mathcal{H})} \underset{n \rightarrow +\infty}{\longrightarrow} 0$. Writing the ISS estimate (\ref{eq: th weak solutions - ISS}) for each classical solution $X_n$ and letting $n \rightarrow +\infty$ shows that (\ref{eq: th weak solutions - ISS}) holds true for $X$ for all $t \in [0,T]$. 

\textit{Step~2: The obtained $X \in \mathcal{C}^0([0,T];\mathcal{H})$ is independent of the chosen approximating sequences of $X_0$, $d$, and $U$.}

For a given $T > 0$, we show that the construction of Step~1 provides a $X \in \mathcal{C}^0([0,T];\mathcal{H})$ that uniquely depends on $(X_0,d,u)$ in the sense that it is independent of the employed approximation sequences. Assume that, following the construction of Step~1, $(X_{1,0,n})_n \in D(\mathcal{A})^{\mathbb{N}}$, $(d_{1,n})_n \in \mathcal{C}^2([0,T];\mathbb{K}^m)^\mathbb{N}$, $(U_{1,n})_n \in \mathcal{C}^1([0,T];\mathcal{H})^\mathbb{N}$ and $(X_{2,0,n})_n \in D(\mathcal{A})^{\mathbb{N}}$, $(d_{2,n})_n \in \mathcal{C}^2([0,T];\mathbb{K}^m)^\mathbb{N}$, $(U_{2,n})_n \in \mathcal{C}^1([0,T];\mathcal{H})^\mathbb{N}$ converge to $X_{0}$, $\left. d \right|_{[0,T]}$, $\left. U \right|_{[0,T]}$, respectively. For any $n\in\mathbb{N}$ and $i \in \{1,2\}$, let $X_{i,n}$ be the unique classical solution associated with $(X_{i,0,n},d_{i,n},U_{i,n})$ over $[0,T]$. We know from Step~1 that $X_{i,n}$ converges to some $X_{i} \in\mathcal{C}^{0}([0,T];\mathcal{H})$ when $n \rightarrow +\infty$. By linearity $X_{1,n} - X_{2,n}$ is the unique classical solution associated with $(X_{1,0,n} - X_{2,0,n},d_{1,n} - d_{2,n},U_{1,n}-U_{2,n})$. Thus the ISS estimate for classical solutions (\ref{eq: th weak solutions - ISS}) yields for all $n\in\mathbb{N}$, 
\begin{align*}
\Vert X_{1,n} - X_{2,n} \Vert_{\mathcal{C}^0([0,T];\mathcal{H})} \leq
& \beta \left( \Vert X_{1,0,n} - X_{2,0,n} \Vert_\mathcal{H} , 0 \right) \\
& + \gamma_1 \left( \Vert d_{1,n} - d_{2,n} \Vert_{\mathcal{C}^0([0,T];\mathbb{K}^m)} \right) \\ 
& + \gamma_2 \left( \Vert U_{1,n} - U_{2,n} \Vert_{\mathcal{C}^0([0,T];\mathcal{H})} \right) .
\end{align*}  
Letting $n \rightarrow +\infty$, it gives $\Vert X_{1} - X_{2} \Vert_{\mathcal{C}^0([0,T];\mathcal{H})} = 0$, i.e., $X_1 = X_2$.

\textit{Step~3: Definition of a weak solution candidate $X \in \mathcal{C}^0(\mathbb{R}_+;\mathcal{H})$.}

Let $0 < T_1 < T_2$ be arbitrarily given and let $X_1 \in \mathcal{C}^0([0,T_1];\mathcal{H})$ and $X_2 \in \mathcal{C}^0([0,T_2];\mathcal{H})$ as provided by the construction of Step~1. It is easy to see that, by restricting the approximation sequences of $\left. d \right|_{[0,T_2]}$ and $\left. U \right|_{[0,T_2]}$ from $[0,T_2]$ to $[0,T_1]$ and by resorting to the uniqueness result of Step~2, that $X_1 = \left. X_2 \right|_{[0,T_1]}$. Therefore, we can define $X \in \mathcal{C}^0(\mathbb{R}_+;\mathcal{H})$ such that for any $T>0$, $\left. X \right|_{[0,T]} \in \mathcal{C}^0([0,T];\mathcal{H})$ is the result of the construction of Step~1. As (\ref{eq: th weak solutions - ISS}) holds true for all $t \in [0,T]$ and for all $T > 0$ with functions $\beta,\gamma_1,\gamma_2$ that are independent of $T$, then (\ref{eq: th weak solutions - ISS}) holds true for the built function $X \in \mathcal{C}^0(\mathbb{R}_+;\mathcal{H})$ for all $t \geq 0$.

\textit{Step~4: The obtained candidate $X \in \mathcal{C}^0(\mathbb{R}_+;\mathcal{H})$ is the unique weak solution associated with $(X_0,d,U)$.}

Let $T > 0$ be arbitrarily given. Let $(X_{0,n})_n \in D(\mathcal{A})^{\mathbb{N}}$, $(d_n)_n \in \mathcal{C}^2([0,T];\mathbb{K}^m)^\mathbb{N}$, and $(U_n)_n \in \mathcal{C}^1([0,T];\mathcal{H})^\mathbb{N}$ be approximating sequences, compliant with the procedure of Step~1, converging to $X_0$, $\left. d \right|_{[0,T]}$, and $\left. U \right|_{[0,T]}$, respectively. Thus, the corresponding sequence of classical solutions $(X_{n})_n$ converges to $\left. X \right|_{[0,T]}$. Based on Corollary~\ref{cor: equiv classical and weak solutions for regular IC and disturbances}, $X_n$ is also a weak solution for all $n\in\mathbb{N}$. Thus, we have for all $z \in \mathcal{C}^0([0,T];D(\mathcal{A}_0^*)) \cap \mathcal{C}^1([0,T];\mathcal{H})$ such that $\mathcal{A}_0^* z \in \mathcal{C}^0([0,T];\mathcal{H})$ and $z(T)=0$, 
\begin{align}
& \int_0^T \left< X_n(t), \mathcal{A}_0^* z(t) + \dfrac{\mathrm{d} z}{\mathrm{d} t}(t) \right>_\mathcal{H} \diff t \nonumber \\
= & - \left< X_{0,n} , z(0) \right>_\mathcal{H} 
+ \int_0^T \left< B d_n(t) , \mathcal{A}_0^* z(t) \right>_\mathcal{H} \diff t \label{eq: weak solution Xn} \\
& - \int_0^T \left< \mathcal{A} B d_n(t) , z(t) \right>_\mathcal{H} \diff t - \int_0^T \left< U_n(t) , z(t) \right>_\mathcal{H} \diff t \nonumber .
\end{align}
From $X_{0,n} \underset{n \rightarrow +\infty}{\longrightarrow} X_0$, we have $\left< X_{0,n} , z(0) \right>_\mathcal{H} \underset{n \rightarrow +\infty}{\longrightarrow} \left< X_{0} , z(0) \right>_\mathcal{H}$. By the Cauchy-Schwarz inequality, we get
\begin{align*}
& \left\vert \int_0^T \left< X_n(t) - X(t), \mathcal{A}_0^* z(t) + \dfrac{\mathrm{d} z}{\mathrm{d} t}(t) \right>_\mathcal{H} \diff t \right\vert \\
& \leq T \left\Vert \mathcal{A}_0^* z + \dfrac{\mathrm{d} z}{\mathrm{d} t} \right\Vert_{\mathcal{C}^0([0,T];\mathcal{H})} \left\Vert X_n - X \right\Vert_{\mathcal{C}^0([0,T];\mathcal{H})} \\
& \underset{n \rightarrow +\infty}{\longrightarrow} 0 .
\end{align*}
Applying a similar procedure to the three integral terms on the right hand side of (\ref{eq: weak solution Xn}), and recalling that operators $B$ and $\mathcal{A} B$ are bounded, one can show their convergence when $n \rightarrow +\infty$. Thus, letting $n \rightarrow +\infty$ in (\ref{eq: weak solution Xn}), we obtain that $X$ satisfies (\ref{eq: def weak solution}) for all $T > 0$ and all test function $z$ over $[0,T]$. Thus, $X$ is the unique weak solution associated with $(X_0,d,U)$. \qed

\end{document}